\documentclass[10pt,sumlimits,intlimits,reqno,oneside]{amsart}
\usepackage{amsmath,amssymb,amsthm,bbm,mathrsfs}
\usepackage[cp1251,cp866]{inputenc}
\usepackage[english,russian]{babel}
\usepackage{delarray}
\inputencoding{cp1251}
\newcounter{nnn}[section]
\theoremstyle{plain}
\newtheorem{Th}{╥хюЁхьр}[section]
\newtheorem{PPP}[Th]{╧Ёхфыюцхэшх}
\newtheorem{Con}[Th]{╤ыхфёЄтшх}

\newtheorem{Nz}[Th]{╟рьхўрэшх}

\newtheorem*{Ex}{╧ЁшьхЁ√}

\begin{document}\thispagestyle{empty}
\numberwithin{equation}{section}
\def\nbb{\refstepcounter{nnn}{\bf\arabic{nnn}.}}
\newcommand{\pr}{
\par\noindent
\nbb{}
} \sloppy \tolerance=10000
\addtolength{\baselineskip}{0.5\baselineskip}
\thispagestyle{empty}
\def\zd{{\mathbb Z}^{d}}
\def\zz{{\mathbb Z}}
\def\rd{{\mathbb R}^{d}}
\def\rr{{\mathbb R}}
\def\NN{{\mathbb N}}
\def\T{\mathbb{T}}
\def\ge{\geqslant}
\def\le{\leqslant}
\def\rho{\varrho}
\def\SS{\mathfrak{S}}
\def\B{\mathfrak{B}}
\def\Ai{\mathrm{Ai}}
\def\A{\mathbb{A}}
\def\H{\mathfrak{H}}
\def\t{\mathfrak{t}}
\def\N{\mathfrak{N}}
\def\X{\mathcal{X}}
\def\Y{\mathcal{Y}}
\def\Z{\mathcal{Z}}
\def\L{\mathcal{L}}
\def\E{\mathcal{E}}
\def\F{\mathcal{F}}
\def\G{\mathcal{G}}
\def\K{\mathcal{K}}
\def\U{\mathcal{U}}
\def\V{\mathbb{V}}
\def\W{\mathbb{W}}
\def\T{\mathbb{T}}
\def\C{\mathbb{C}}
\def\supp{\mathrm{supp}}
\def\No{╣}

\title[└ёшьяЄюЄшър фшёъЁхЄэюую ёяхъЄЁр]
{└ёшьяЄюЄшър ш юЎхэъш фшёъЁхЄэюую ёяхъЄЁр \\  юяхЁрЄюЁр ╪ЁхфшэухЁр
эр фшёъЁхЄэюь яхЁшюфшўхёъюь уЁрЇх}

\date{\today}

\author{┼. ╦. ╩юЁюЄ хт$^{*}$, ┬. └. ╤ыює∙$^{**}$}
\thanks{{\it ╩ы■ўхт√х ёыютр:\/} фшёъЁхЄэ√щ юяхЁрЄюЁ ╪ЁхфшэухЁр, шэЄхуЁры№э√х юяхЁрЄюЁ√,
 юЎхэъш ёшэуєы Ёэ√ї ўшёхы, ъырёё√ ъюьяръЄэ√ї юяхЁрЄюЁют.}
\thanks{$^{*}$ ╨рсюЄр ртЄюЁр яюффхЁцрэр уЁрэЄюь ╨═╘ 18-11-00032}
\thanks{$**$ ╨рсюЄр ртЄюЁр яюффхЁцрэр уЁрэЄюь ╨╘╘╚ 17-01-00668}
\address{╤.-╧хЄхЁсєЁуёъшщ\\ уюёєфрЁёЄтхээ√щ єэштхЁёшЄхЄ\\ ърЇхфЁр ьрЄхьрЄшўхёъюую рэрышчр $^*$,
ърЇхфЁр т√ё°хщ ьрЄхьрЄшъш  ш ьрЄхьрЄшўхёъющ Їшчшъш $^{**}$\\
╙эштхЁёшЄхЄёър  эрсхЁхцэр  7/9\\ 198034, ╤рэъЄ-╧х\-ЄхЁ\-сєЁу\\
╨юёёш } \email{e.korotyaev@spbu.ru, korotyaev@gmail.com,
v.slouzh@spbu.ru, vsloushch@list.ru}

\begin{abstract}
╨рёёьрЄЁштрхЄё  яхЁшюфшўхёъшщ юяхЁрЄюЁ ╪ЁхфшэухЁр $H$ эр фшёъЁхЄэюь
яхЁшюфшўхёъюь уЁрЇх.  ╧юыєўхэ√ юЎхэъш фшёъЁхЄэюую ёяхъЄЁр
тючьє∙хээюую юяхЁрЄюЁр $H_{\pm}(t)=H\pm tV$, $t>0$, уфх  $V\ge 0$
єс√тр■∙шщ яюЄхэЎшры. ┬ ёыєўрх яюЄхэЎшрыр, шьх■∙хую ёЄхяхээє■
рёшьяЄюЄшъє эр схёъюэхўэюёЄш
эрщфхэр
рёшь\-я\-Єю\-Єшър фшёъЁхЄэюую ёяхъЄЁр юяхЁрЄюЁр $H_{\pm}(t)$ яЁш
сюы№°ющ ъюэёЄрэЄх ёт чш.

\end{abstract}
\maketitle \thispagestyle{empty}
\setcounter{section}{-1}

\section{┬тхфхэшх}\setcounter{nnn}{0}

\noindent\nbb\noindent{} ╬яхЁрЄюЁ√ ╦ряырёр ш ╪ЁхфшэухЁр эр уЁрЇрї
шьх■Є ьэюуюўшёыхээ√х яЁшыюцхэш  т Їшчшъх ш їшьшш, ёь., эряЁшьхЁ,
юсчюЁ \cite{CGPNG09} ш ёё√ыъш т эхь (фтюх ртЄюЁют юсчюЁр
═ютюёхыютр ш ├хщьр яюыєўшыш эюсхыхтёъє■ яЁхьш■ чр юЄъЁ√Єшх
уЁрЇхэр).
═р фшёъЁхЄэюь уЁрЇх юяхЁрЄюЁ ╦ряырёр фхщёЄтєхЄ т яЁюёЄЁрэёЄтх
ЇєэъЎшщ, юяЁхфхыхээ√ї эр ьэюцхёЄтх тхЁ°шэ уЁрЇр. ╟фхё№ уыртэє■
Ёюы№ шуЁр■Є тхЁ°шэ√ уЁрЇр, р ЁхсЁр яюърч√тр■Є тчршьюфхщёЄтшх ьхцфє
тхЁ°шэрьш уЁрЇр. ┬ юёэютэюь шчєўр■Є фтр тшфр {\it фшёъЁхЄэюую\/}
юяхЁрЄюЁр ╦ряырёр: ъюьсшэрЄюЁэ√щ юяхЁрЄюЁ ╦ряырёр
\cite{KorSlo.KS2014} ш эюЁьшЁютрээ√щ юяхЁрЄюЁ ╦ряырёр
\cite{KorSlo.KS2018}; шчєўр■Є ш фЁєушх тшф√ юяхЁрЄюЁют ╦ряырёр эр
фшёъЁхЄэюь уЁрЇх (яю яютюфє Ёрчышўэ√ї юяхЁхфхыхэшщ юяхЁрЄюЁр
╦ряырёр эр уЁрЇх ёь. Єръцх \cite{Ch97}, \cite{S13}). ╤яхъЄЁ
фшёъЁхЄэюую юяхЁрЄюЁр ╦ряырёр эр яхЁшюфшўхёъюь уЁрЇх ёюёЄюшЄ шч
ъюэхўэюую ўшёыр чюэ, Ёрчфхыхээ√ї ыръєэрьш.


┴юы№°шэёЄтю Ёхчєы№ЄрЄют фы  фшёъЁхЄэ√ї юяхЁрЄюЁют ╪Ё╕фшэухЁр
$H=\Delta +V$
фы  ЁртэюьхЁэю єс√тр■∙шї яюЄхэЎшрыют $V$ яюыєўхэю эр ьэюуюьхЁэющ
ъєсшўхёъющ Ёх°хЄъх $\zd$, уфх $\Delta$ --- юяхЁрЄюЁ ╦ряырёр эр
$\ell_2(\zd)$. ┴єЄх фх ╠юэтхы№ ш ╤рїсрэш \cite{BS99} шчєўрыш
ёыєўрщ ёЄхяхээюую єс√трэш  $V$ эр схёъюэхўэюёЄш (яЁш яюърчрЄхых
єс√трэш  сюы№°х $2$).
╬эш яюърчрыш, шёяюы№чє  ьхЄюф ╠єЁЁр, ўЄю тюыэют√х юяхЁрЄюЁ√
$$
W_\pm = s-\lim_{} e^{itH}e^{-it\Delta},\ \ \ \ \ t\to \pm\infty,
$$
ёє∙хёЄтє■Є ш яюыэ√, Є.х., Ran $W_\pm={\mathscr H}_{ac}(H)$. ┴юыхх
Єюую, юэш фюърчрыш юЄёєЄёЄтшх ёшэуєы Ёэюую эхяЁхЁ√тэюую ёяхъЄЁр ш
ыюъры№эє■ ъюэхўэюёЄ№ ёюсёЄтхээ√ї чэрўхэшщ $H$ тэх ъЁрхт
эхяЁхЁ√тэюую ёяхъЄЁр. ┬тхфхь юяхЁрЄюЁ-ЇєэъЎш■
$$
F(z)=|V|^{\frac{1}{2}}(\Delta-z)^{-1}|V|^{\frac{1}{2}}, \ \ z\in
\Lambda=\C\setminus \sigma(\Delta).
$$
┬ ЁрсюЄрї \cite{BS99}, \cite{SW01}, яюърчрэю, ўЄю юяхЁрЄюЁ-ЇєэъЎш 
$ F(z)$ рэрышЄшўхёър  т $\Lambda$  ш эхяЁхЁ√тэр тяыюЄ№ фю уЁрэшЎ√
чр шёъы■ўхэшхь ъюэхўэюую ўшёыр Єюўхъ эр ёяхъЄЁх $\sigma(\Delta)$.
╚ёючръш, ╩юЁюЄ хт \cite{IK12}  яюърчрыш, ўЄю Єръшї Єюўхъ эхЄ яЁш
$d\ge 3$. ┴юыхх Єюую юэш Ёх°шыш юсЁрЄэє■ чрфрўє ю тюёёЄрэютыхэшш
яюЄхэЎшрыр яю ьрЄЁшЎх Ёрёёх эш  яЁш тёхї ¤эхЁуш ї, р Єръцх
яюыєўшыш ЇюЁьєы√ ёыхфют. ╚чюёръш, ╠юЁшюъш фюърчрыш, ўЄю Єюўхўэ√щ
ёяхъЄЁ юяхЁрЄюЁют ╪Ё╕фшэухЁр эр эхяЁхЁ√тэюь ёяхъЄЁх юЄёєЄёЄтєхЄ т
ёыєўрх яюЄхэЎшрыр ё ъюэхў√ь эюёшЄхыхь. ╤ыєўрщ яюЄхэЎшрыр $V\in
\ell_p(\mathbb{Z}^d)$ шчєўрыё  т ЁрсюЄх \cite{KM17}. ╟фхё№
ЁрёёьюЄЁхэ юяхЁрЄюЁ ╦ряырёр эр Ёх°хЄъх $\mathbb{Z}^d$ ш яюыєўхэ√
юЎхэъш уЁєяя√ $e^{it\Delta}$ ш Ёхчюы№тхэЄ√ $(\Delta-z)^{-1}$,
фхщёЄтє■∙шї т яЁюёЄЁрэёЄтрї $\ell_q(\mathbb{Z}^d)$. ▌Єш Ёхчєы№ЄрЄ√
яЁшьхэ ■Єё  ъ шчєўхэш■ юяхЁрЄюЁр ╪Ё╕фшэухЁр т
$\ell_2(\mathbb{Z}^d)$ ё яюЄхэЎшрыюь шч $ \ell_p(\mathbb{Z}^d)$
яЁш яюфїюф ∙хь $p\ge 1$.
 ╤юсёЄтхээ√х
чэрўхэш  ш ЇюЁьєы√ ёыхфют фы  ъюьяыхъёэюую яюЄхэЎшрыр
ЁрёёьрЄЁштрышё№ т \cite{KL18}, \cite{K17}. ╧рЁЁр ш ╨шўрЁф
\cite{PR18} ЁрёяЁюёЄЁрэшыш Ёхчєы№ЄрЄ√ ┴єЄх фх ╠юэтхы , ╤рїсрэш
\cite{BS99} эр юс∙шх яхЁшюфшўхёъшх  уЁрЇ√.


═рёЄю ∙р  ЁрсюЄр яюёт ∙хэр  юяшёрэш■ фшёъЁхЄэюую ёяхъЄЁр
$\zd$-яхЁшюфшўхёъюую юяхЁрЄюЁр ╪ЁхфшэухЁр, тючьє∙хээюую єс√тр■∙шь
яюЄхэЎшрыюь, эр  фшёъЁхЄэюь $\zd$-яхЁшюфшўхёъюь уЁрЇх.
╧юфЁюсэхх, т ърўхёЄтх эхтючьє∙хээюую юяхЁрЄюЁр ЁрёёьюЄЁшь
фшёъЁхЄэ√щ яхЁшюфшўхёъшщ юяхЁрЄюЁ ╪ЁхфшэухЁр $H=\Delta+Q$ эр
ёт чэюь ыюъры№эю-ъюэхўэюь $\zd$-яхЁшюфшўхёъюь уЁрЇх $G$, тыюцхээюь
т $\rd$; чфхё№ $\Delta$
--- фшёъЁхЄэ√щ (ъюьсшэрЄюЁэ√щ) юяхЁрЄюЁ ╦ряырёр эр $G$, $Q$
--- тх∙хёЄтхээ√щ юуЁрэшўхээ√щ $\zd$-яхЁшюфшўхёъшщ яюЄхэЎшры эр
$G$. ╤яхъЄЁ юяхЁрЄюЁр $H$ ёюёЄюшЄ шч ъюэхўэюую ўшёыр чюэ,
Ёрчфхыхээ√ї ыръєэрьш; эхъюЄюЁ√х чюэ√ (эх тёх ёь.
\cite{KorSlo.KS2016}) ьюуєЄ т√ЁюцфрЄ№ё  т Єюўъє (ёюсёЄтхээюх
чэрўхэшх юяхЁрЄюЁр $H$ схёъюэхўэющ ъЁрЄэюёЄш).

╨рёёьюЄЁшь тючьє∙хээ√х  юяхЁрЄюЁ√ ╪ЁхфшэухЁр $H_{\pm}(t):=H\pm
tV$, $t>0$, уфх $V$ --- эхюЄЁшЎрЄхы№э√щ єс√тр■∙шщ эр схёъюэхўэюёЄш
яюЄхэЎшры.
═рё шэЄхЁхёєхЄ ёяхъЄЁ юяхЁрЄюЁют $H_{\pm}(t)$, тючэшър■∙шщ т
ыръєэрї ёяхъЄЁр юяхЁрЄюЁр $H$. ╧єёЄ№ $(\Lambda_{+},\Lambda_{-})$
 ты хЄё  ыръєэющ (тючьюцэю яюыєсхёъюэхўэющ) т ёяхъЄЁх юяхЁрЄюЁр
$H$. ╧юёъюы№ъє яюЄхэЎшры $V$ єс√трхЄ эр схёъюэхўэюёЄш, юяхЁрЄюЁ
єьэюцхэш  эр $V$ ъюьяръЄхэ, р яюЄюьє ёяхъЄЁ юяхЁрЄюЁр
$H_{\pm}(t)$, $t>0$, т ыръєэх $(\Lambda_{+},\Lambda_{-})$
фшёъЁхЄхэ. ╤юсёЄтхээ√х чэрўхэш  юяхЁрЄюЁют $H_{\pm}(t)$ ьюэюЄюээю
фтшцєЄё  ё ЁюёЄюь $t$. ╬ёэютэ√ь юс·хъЄюь шёёыхфютрэш   ты ■Єё 
ёўшЄр■∙шх ЇєэъЎшш $N_{_{\pm}}(\lambda,\tau)$,
$\lambda\in(\Lambda_{+},\Lambda_{-})$, $\tau>0$, Ёртэ√х ўшёыє
ёюсёЄтхээ√ї чэрўхэшщ юяхЁрЄюЁют $H_{\pm}(t)$, яЁю°хф°шї ўхЁхч
ЇшъёшЁютрээє■  Єюўъє $\lambda$ яЁш єтхышўхэшш $t$ юЄ $0$ фю
$\tau$. ╤ўшЄр■∙шх ЇєэъЎшш эр ъЁр■ ыръєэ√, юяЁхфхыхээ√х ъръ яЁхфхы
$$
N_{\pm}(\Lambda_{\pm},\tau):=\lim\limits_{\lambda\to\Lambda_{\pm}\pm
0}N_{\pm}(\lambda,\tau),
$$
ьюуєЄ шэЄхЁяЁхЄшЁютрЄ№ё , ъръ ўшёыю ёюсёЄтхээ√ї чэрўхэшщ,
<<Ёюфшт°шїё >> т Єюўъх $\Lambda_{\pm}$. ╨рчєьххЄё , хёыш
$(\Lambda_{+},+\infty)$ --- яюыєсхёъюэхўэр  ыръєэр т ёяхъЄЁх
юяхЁрЄюЁр $H$, Єю ўшёыю $N_{+}(\lambda,\tau)$,
$\lambda\in[\Lambda_{+},+\infty)$, хёЄ№ ёєььрЁэр  ъЁрЄюёЄ№
ёюсёЄтхээ√ї чэрўхэшщ юяхЁрЄюЁр $H_{+}(\tau)$ эр шэЄхЁтрых
$(\lambda,+\infty)$. └эрыюушўэю шэЄхЁяЁхЄшЁєхЄё  тхышўшэр
$N_{-}(\lambda,\tau)$, $\lambda\in(-\infty,\Lambda_{-}]$, хёыш
$(-\infty,\Lambda_{-})$ --- яюыєсхёъюэхўэр  ыръєэр т ёяхъЄЁх
$\sigma(H)$. ╤єььрЁэр  ъЁрЄэюёЄ№ юяхЁрЄюЁр $H_{i}(\tau)$ тю {\it
тэєЄЁхээхщ\/} ыръєэх $(\Lambda_{+},\Lambda_{-})$ эх яЁхтюёїюфшЄ
тхышўшэ√ $N_{i}(\Lambda_{i},\tau)$, $i=\pm$.


╬ёэютэ√х Ёхчєы№ЄрЄ√ ЁрсюЄ√  ёюёЄю Є т ёыхфє■∙хь.

\noindent {\it {\rm 1)} ┬  ЄхюЁхьх
{\rm\ref{KorSloTheoremAsymptotic}} яюърчрэю, ўЄю яЁш єёыютшш
\begin{equation}\label{KorSlo.0.1}
\textstyle
0\le V(x)\sim\vartheta(\frac{x}{|x|})|x|^{-d/p},\ \ |x|\to\infty,\
\ p>0,
\end{equation}
ёўшЄр■∙шх ЇєэъЎшш шьх■Є ёЄхяхээє■ рёшьяЄюЄшъє яю сюы№°ющ ъюэёЄрэЄх
ёт чш яЁш тёхї ЇшъёшЁютрээ√ї $\lambda\in(\Lambda_{+},\Lambda_{-})$
\begin{equation}\label{KorSlo.0.2}
N_{\pm}(\lambda,\tau)=\tau^{p}(\Gamma_{p}^{\pm}(\lambda)+o(1)),\ \ \
\ \tau\to+\infty,
\end{equation}
уфх $\Gamma_{p}^{\pm}(\lambda)=\Gamma_{p}^{\pm}(\lambda,H,V)$.
╧Ёш
юяЁхфхыхээ√ї єёыютш ї рёшьяЄюЄшър {\rm(\ref{KorSlo.0.2})} шьххЄ
ьхёЄю ш эр ъЁр ї ыръєэ√.
\par\noindent
{\rm 2)} ┬ ЄхюЁхьрї {\rm\ref{KorSloTheoremRLC}},
{\rm\ref{KorSloTheoremWRLC}} яЁштхфхэ√ єёыютш  ъюэхўэюёЄш ёяхъЄЁр
юяхЁрЄюЁют $H_{\pm}(t)$, $t>0$, т ёяхъЄЁры№э√ї ыръєэрї юяхЁрЄюЁр
$H$ ш фрэ√ юЎхэъш ёєььрЁэющ ъЁрЄэюёЄш ¤Єюую ёяхъЄЁр.}
\pr ─шёъЁхЄэ√щ ёяхъЄЁ яхЁшюфшўхёъюую юяхЁрЄюЁр ╪ЁхфшэухЁр,
тючьє∙хээюую єс√тр■∙шь яюЄхэЎшрыюь, шчєўрыё  тю ьэюушї ЁрсюЄрї
(ёь., эряЁшьхЁ, юсчюЁ√ \cite{KorSlo.Birman96},
\cite{KorSlo.RozenblumSolomyak2008}, р Єръцх \cite{KorSlo.Bach}).
╬ёЄрэютшьё , яЁхцфх тёхую, эр {\it эхяЁхЁ√тэюь\/} ёыєўрх, Є.х. эр
ёыєўрх тючьє∙хэш  юяхЁрЄюЁр ╪ЁхфшэухЁр $H=-\Delta+Q(x)$ т
$L_{2}(\rd)$ єс√тр■∙шь яюЄхэЎшрыюь $V(x)$. ┬ эхяЁхЁ√тэюь ёыєўрх
їрЁръЄхЁ Ёхчєы№ЄрЄют ш шёяюы№чєхьр  Єхїэшър ёє∙хёЄтхээю чртшё Є юЄ
чэрър тючьє∙хэш .

─ы  {\it юЄЁшЎрЄхы№э√ї\/} тючьє∙хэшщ юяхЁрЄюЁр ╦ряырёр
╨ю\-чхэ\-сы■\-ьюь  \cite{KorSlo.Rozenblum72},
\cite{KorSlo.Rozenblum76}, ╦ш\-сюь \cite{KorSlo.Lieb76} ш
╓тш\-ъх\-ыхь \cite{KorSlo.Cwikel77} с√ы яюыєўхэ ёыхфє■∙шщ
Ёхчєы№ЄрЄ.
\par\noindent{\it ╧єёЄ№ $H=-\Delta$ т $L_{2}(\rd)$, ш яєёЄ№ тючьє∙хэшх $V$ єфютыхЄтюЁ хЄ єёыютш■
\begin{equation}\label{KorSlo.0.3}
0\le V\in L_{d/2}(\rd),\ \ d\ge 3.
\end{equation}
╥юуфр фы  ы■сюую $\lambda\le 0$  ёяЁртхфышт√ юЎхэър ш рёшьяЄюЄшър
\begin{equation}\label{KorSlo.0.4}
N_{-}(\lambda,\tau)\le C(d)\tau^{d/2}\|V\|_{L_{\frac d2}}^{\frac
d2},\ \ \tau>0,\ \
\end{equation}
\begin{equation}\label{KorSlo.0.5}
N_{-}(\lambda,\tau)=\omega_{d}\tau^{d/2}(\|V\|_{L_{d/2}}^{d/2}+o(1)),\
\ \tau\to+\infty,\
\end{equation}
уфх $\omega_{d}$ --- юс·хь хфшэшўэюую °рЁр т $\rd$.\/}

 ┬ ЁрсюЄрї ╒хьяхы  \cite{KorSlo.Hempl89} ш ┴шЁьрэр \cite{KorSlo.Birman91} рёшьяЄюЄшър
 (\ref{KorSlo.0.5}) с√ыр фюърчрэр фы  яхЁшюфшўхёъюую юяхЁрЄюЁр ╪ЁхфшэухЁр
 $H=-\Delta+Q(x)$  фы 
 яЁюшчтюы№эюую $\lambda\in(\Lambda_{+},\Lambda_{-})$ яЁш єёыютшш (\ref{KorSlo.0.3}).
  ╧Ёш ¤Єюь т√ ёэшыюё№, ўЄю рёшьяЄюЄрЄшўхёъшщ ъю¤ЇЇшЎшхэЄ эх чртшёшЄ юЄ Єюўъш эрсы■фхэш 
  $\lambda\in(\Lambda_{+},\Lambda_{-})$ ш яюЄхэЎшрыр
 $Q(x)$. ╤ыєўрщ рёшьяЄюЄшъш эр ъЁр■ ёяхъЄЁры№эющ ыръєэ√ яЁюшчтюы№эюую яхЁшюфшўхёъюую
 юяхЁрЄюЁр ╪ЁхфшэухЁр юсёєцфрыё  т ЁрсюЄх
 \cite{KorSlo.Birman95}.

 ┴юы№°юх ўшёыю ЁрсюЄ (яюфЁюсэхх ёь., эряЁшьхЁ, \cite{KorSlo.Birman96}, \cite{KorSlo.RozenblumSolomyak2008}, \cite{Sha}) яюёт ∙хэ√ юсёєцфхэш■ {\it эхЁхуєы Ёэ√ї
 тючьє∙хэшщ\/} фы  ъюЄюЁ√ї єёыютшх (\ref{KorSlo.0.3}) эх т√яюыэхэю (т Є.ў. яЁш $d=1,2$).

 ═хюЄЁшЎрЄхы№э√х тючьє∙хэш  шчєўрышё№ т ЁрсюЄрї └ырь√, ─хщЇЄр, ╒хьяы  \cite{KorSlo.ADH}, ┴шЁьрэр
 \cite{KorSlo.Birman91}, \cite{KorSlo.BirmanSloushch}, \cite{KorSlo.BirmanM94}, ш ╤ыює∙р
\cite{KorSlo.SloushchDiscreteSpectrum}. ╧Ёш ¤Єюь фы  яхЁшюфшўхёъюую
юяхЁрЄюЁр ╪ЁхфшэухЁр $H=-\Delta+Q(x)$ с√ы яюыєўхэ ёыхфє■∙шщ
Ёхчєы№ЄрЄ.
\par\noindent
{\it ┼ёыш тючьє∙хэшх $V$ єфютыхЄтюЁ хЄ єёыютш■
{\rm(\ref{KorSlo.0.1})}, Єю ёўшЄр■∙р  ЇєэъЎш  $N_{+}(\lambda,\tau)$
шьххЄ ёЄхяхээє■ рёшьяЄюЄшъє {\rm(\ref{KorSlo.0.2}$+$)} яю сюы№°ющ
ъюэёЄрэЄх ёт чш яЁш тёхї $\lambda\in(\Lambda_{+},\Lambda_{-})$. ╧Ёш
юяЁхфхыхээ√ї єёыютш ї рёшьяЄюЄшър {\rm(\ref{KorSlo.0.2}$+$)}
ёяЁртхфыштр ш яЁш $\lambda=\Lambda_{+}$.\/}

╬Ўхэър тхышўшэ√ $N_{+}(\lambda,\tau)$ юсёєцфрырё№ т ЁрсюЄх
\cite{KorSlo.BirmanM94}.

┬ючьє∙хэш■ яхЁшюфшўхёъюую {\it фшёъЁхЄэюую\/} юяхЁрЄюЁр ╪ЁхфшэухЁр
єс√тр■∙шь яюЄхэЎшрыюь яюёт ∙хэю сюы№°юх ъюышўхёЄтю ЁрсюЄ (ёь.
юсчюЁ Ёхчєы№ЄрЄют т \cite{KorSlo.RozenblumSolomyak2008} ш
\cite{KorSlo.Bach}). ╥ръ т ЁрсюЄрї ╦хтшэр, ╨ючхэсы■ьр ш ╤юыюь ър
\cite{KorSlo.RozenblumSolomyak2009}, \cite{KorSlo.LS97} ш
\cite{KorSlo.RozenblumSolomyak2008} шчєўрыюё№ тючьє∙хэшх
фшёъЁхЄэюую юяхЁрЄюЁр ╦ряырёр єс√тр■∙шь яюЄхэЎшрыюь; с√ыш яюыєўхэ√
юЎхэъш тхышўшэ√ $N_{-}(0,\tau)$, т Є.ў. юЎхэър (\ref{KorSlo.0.4})
яЁш єёыютшш $0\le V\in L_{d/2}(G)$, $d\ge 3$; юЎхэъш тхышўшэ√
$N_{-}(0,\tau)$ яЁш $d=2$ шчєўрышё№ т ЁрсюЄрї \cite{MolVain} ш
\cite{RozSol2014}.

┬ ЁрсюЄрї \cite{KorSlo.RozenblumSolomyak2010} ш \cite{KorSlo.Bach}
с√ыш єёЄрэютыхэ√ юЎхэъш тхышўшэ√ $N_{-}(0,\tau)$ т ёыєўрх
тючьє∙хэш  фшёъЁхЄэюую юяхЁрЄюЁр ╪ЁхфшэухЁр єс√тр■∙шь яюЄхэЎшрыюь.
└ёшьяЄюЄшър тхышўшэ√ $N_{-}(0,\tau)$ яю сюы№°ющ ъюэёЄрэЄх ёт чш
$\tau$ юсёєцфрырё№ т \cite{KorSlo.Bach} фы  ёыєўр  тючьє∙р■∙хую
яюЄхэЎшрыр єс√тр■∙хую эр схёъюэхўэюёЄш ётхЁїёЄхяхээ√ь юсЁрчюь.
 ╒р ёш, ╒шуєўш,
═юьєЁр ш ╬уєЁшёє \cite{HHNO16} т√ўшёышыш ўшёыю фшёъЁхЄэ√ї
ёюсёЄтхээ√ї чэрўхэшщ фы  юяхЁрЄюЁр ╦ряырёр, тючьє∙хээюую
яюЄхэЎшрыюь ё ъюэхўэ√ь эюёшЄхыхь. ┬√сЁрээ√щ эрьш рёшьяЄюЄшўхёъшщ
Ёхцшь (\ref{KorSlo.0.1}), (\ref{KorSlo.0.2}) яю-тшфшьюьє Ёрэхх эх
ЁрёёьрЄЁштрыё .
\pr
╧ю їрЁръЄхЁє яюыєўхээ√ї Ёхчєы№ЄрЄют ш шёяюы№чютрээющ Єхїэшъх эр°р
ЁрсюЄр сышчър ъ ЁрсюЄрь ┴шЁьрэр ш ╤ыює∙р
\cite{KorSlo.BirmanSloushch},
\cite{KorSlo.SloushchDiscreteSpectrum}, уфх шчєўрыё  фшёъЁхЄэ√щ
ёяхъЄЁ юяхЁрЄюЁр ╪ЁхфшэухЁр т $L_{2}(\rd)$, тючьє∙хээюую {\it
эхюЄЁшЎрЄхы№э√ь\/} яюЄхэЎшрыюь єс√тр■∙шь эр схёъюэхўэюёЄш ёЄхяхээ√ь
юсЁрчюь.
╩ръ ш т ЁрсюЄрї \cite{KorSlo.BirmanSloushch},
\cite{KorSlo.SloushchDiscreteSpectrum} ь√ ётюфшь (ёь. яЁхфыюцхэшх
\ref{KorSloPropositionBirmanShwinger}) шчєўхэшх рёшьяЄюЄшъш
тхышўшэ $N_{\pm}(\lambda,\tau)$ яю сюы№°ющ ъюэёЄрэЄх ёт чш $\tau$
ъ шёёыхфютрэш■ юЎхэюъ ш рёшьяЄюЄшъ ёшэуєы Ёэ√ї ўшёхы яюфїюф ∙шї
<<фшёъЁхЄэ√ї>> ╧─╬ (ёь. \S 2) юЄЁшЎрЄхы№эюую яюЁ фър. ╬ёэютэюх
юЄышўшх эрёЄю ∙хщ ЁрсюЄ√ юЄ \cite{KorSlo.BirmanSloushch} ш
\cite{KorSlo.SloushchDiscreteSpectrum} ёюёЄюшЄ т ёыхфє■∙хь: яЁш
Ёрчыюцхэшш т яЁ ьющ шэЄхуЁры яхЁшюфшўхёъюую фшёъЁхЄэюую юяхЁрЄюЁр
╪ЁхфшэухЁр ёыюш яюыєўр■Єё  ъюэхўэюьхЁэ√ьш; яЁш Ёрчыюцхэшш т яЁ ьющ
шэЄхуЁры яхЁшюфшўхёъюую фшЇЇхЁхэЎшры№эюую юяхЁрЄюЁр ёыюш
схёъюэхўэюьхЁэ√. ╧юёыхфэхх юсёЄю Єхы№ёЄтю фхырхЄ чрфрўє,
ЁрёёьюЄЁхээє■ т \cite{KorSlo.BirmanSloushch},
\cite{KorSlo.SloushchDiscreteSpectrum}, чэрўшЄхы№эю сюыхх ёыюцэющ,
ўхь чрфрўр, шёёыхфєхьр  т эрёЄю ∙хщ ЁрсюЄх.

╧Ёш яюыєўхэшш юЎхэюъ $s$-ўшёхы <<фшёъЁхЄэ√ї>> ╧─╬ юЄЁшЎрЄхы№эюую
яюЁ фър ь√ яюы№чєхьё  шчтхёЄэ√ьш Ёхчєы№ЄрЄрьш ┴шЁьрэр, ╤юыюь ър
\cite{KorSlo.BirmanSolomyak90}.
╧Ёш шёёыхфютрэшш рёшьяЄюЄшъш $s$-ўшёхы <<фшёъЁхЄэ√ї>> ╧─╬
юЄЁшЎрЄхы№эюую яюЁ фър ь√ ёыхфєхь ЁрсюЄрь ┴шЁьрэр, ╤юыюь ър
\cite{KorSlo.BirmanSolomyak77}, \cite{KorSlo.BirmanSolomyak79} ш
ётюфшь фхыю ъ ёыєўр■ <<эхяЁхЁ√тэюую>> ╧─╬. ╥ръцх юЄьхЄшь, ўЄю т
эр°хщ ЁрсюЄх ёє∙хёЄтхээю шёяюы№чє■ё  Ёхчєы№ЄрЄ√ шч
\cite{KorSlo.KS2014} ю яхЁшюфшўхёъюь юяхЁрЄюЁх  ╪ЁхфшэухЁр эр
яхЁшюфшўхёъюь уЁрЇх.
%
\pr
╙ёыютшьё  ю эхъюЄюЁ√ї юсючэрўхэш ї. ╧єёЄ№ $(\X,d\rho)$ ---
шчьхЁшьюх яЁюёЄЁрэёЄтю; фы  шчьхЁшьющ эр $\X$ ЇєэъЎшш $f$
ўхЁхч $[f(x)]$ (р шэюуфр ўхЁхч $f$) юсючэрўрхЄё  юяхЁрЄюЁ
єьэюцхэш  эр $f$ т яЁюёЄЁрэёЄтх $L_{2}(\X,d\rho)$; ъЁюьх Єюую, ь√
шёяюы№чєхь юсючэрўхэш 
\begin{equation*}
f_{\pm}(x)=(|f(x)|\pm f(x))/2,\ \ x\in\X;\ \ f_{\pm}^{-p}(x)=
\left\{
\begin{array}{l}
(f_{\pm}(x))^{-p},\ \ f_{\pm}(x)>0,\\
0,\ \ \ \ \ \ \ \ \ \ \ \ \ \ f_{\pm}(x)=0,
\end{array}
\right. \ \ x\in\X,\ \ p>0.
\end{equation*}
 ╫хЁхч $\mathbf{1}_{\U}$ юсючэрўшь їрЁръЄхЁшёЄшўхёъє■ ЇєэъЎш■
шчьхЁшьюую ьэюцхёЄтр $\U\subset\X$.

═юЁьр т (ътрчш)эюЁьшЁютрээюь яЁюёЄЁрэёЄтх $X$ юсючэрўрхЄё  ўхЁхч
$\|\cdot\|_{X}$. ─ы  ёхярЁрсхы№эюую шчьхЁшьюую яЁюёЄЁрэёЄтр ё
$\sigma$-ъюэхўэющ ьхЁющ $(\Z,d\mu)$ ўхЁхч $L_{p}(\Z,d\mu)$
юсючэрўрхЄё  ёЄрэфрЁЄэ√щ $L_{p}$-ъырёё; фы  яЁюшчтюы№эюую
ьэюцхёЄтр $K$ ўхЁхч $B(K)$ юсючэрўшь ьэюцхёЄтю юуЁрэшўхээ√ї эр $K$
ЇєэъЎшщ; $\|f\|_{B}:=\sup\limits_{k\in K}|f(k)|$, $f\in B(K)$.

─ы  яыюЄ\-эю юя\-Ёх\-фх\-ыхэ\-эюую чрьъэєЄюую юяхЁрЄюЁр $A$ т
ушы№\-схЁ\-Єю\-тюь яЁюёЄЁрэёЄтх ўхЁхч $\sigma(A)$ юсючэрўрхЄё 
ёяхъЄЁ, ўхЁхч $A^{*}$ --- ёюяЁ цхээ√щ юяхЁрЄюЁ. ┼ёыш $A$
ёрьюёюяЁ цхэ, Єю $E_{A}(\cdot)$ --- ёяхъЄЁры№эр  ьхЁр юяхЁрЄюЁр
$A$. ─ы  сюЁхыхтёъюую ьэюцхёЄтр $\delta\subset\mathbb{R}$ яюырурхь
$\pi_{A}(\delta):=\mathrm{rank}E_{A}(\delta)$. ─рыхх,
$n_{\pm}(s,A):=\pi_{\pm A}(s,+\infty)$, $s>0$, --- ЇєэъЎшш
ЁрёяЁхфхыхэш  яюыюцшЄхы№эюую ш юЄЁшЎрЄхы№эюую ёяхъЄЁр юяхЁрЄюЁр
$A$.

╩ЁрЄъю юяш°хь яырэ ЁрсюЄ√, ъюЄюЁр  ёюёЄюшЄ шч ттхфхэш  ш ЄЁхї
ярЁруЁрЇют.

┬ \S 1 яЁштхфхэ√ яюёЄрэютър чрфрўш ш юяшёрэшх юёэютэ√ї
Ёхчєы№ЄрЄют.

┬ \S 2 фрэ√ эхюсїюфшь√х ётхфхэш  юс юЎхэърї ш рёшьяЄюЄшърї $s$-ўшёхы
<<фшёъЁхЄэ√ї>> ╧─╬.

┬ \S 3 яЁштхфхэ√ фюърчрЄхы№ёЄтр юёэютэ√ї ЄхюЁхь
\ref{KorSloTheoremRLC}--\ref{KorSloTheoremAsymptotic}.

└тЄюЁ√ т√Ёрцр■Є ётю■ яЁшчэрЄхы№эюёЄ№ └.┬.┴рфрэшэє, └.╚.═рчрЁютє,
═.─.╘шыюэютє  ш ═.▐.╤рсєЁютющ чр яюыхчэ√х юсёєцфхэш .

\section{╬ёэютэющ Ёхчєы№ЄрЄ}\setcounter{nnn}{0}

\noindent\nbb{\bf{}  ═хтючьє∙хээ√щ юяхЁрЄюЁ.} ╧єёЄ№ т
$\mathbb{R}^{d}$ чрфрэ ёт чэ√щ уЁрЇ $G=(X,\mathcal{E})$, тючьюцэю
шьх■∙шщ яхЄыш шыш ъЁрЄэ√х ЁхсЁр. ╟фхё№ $X$ --- ьэюцхёЄтю тхЁ°шэ
уЁрЇр $G$, $\mathcal{E}$ --- ьэюцхёЄтю эхюЁшхэЄшЁютрээ√ї ЁхсхЁ
уЁрЇр $G$. ╨хсЁю $\mathbf{e}$, ёюхфш ■∙хх тхЁ°шэ√ $x,y\in X$,
эрчютхь ЁхсЁюь {\it шэЎшфхэЄэ√ь\/} тхЁ°шэрь $x$ ш $y$. ┬хЁ°шэ√
$x,y\in X$, ёюхфшэхээ√х ЁхсЁюь $\mathbf{e}\in\mathcal{E}$, эрчютхь
{\it ёьхцэ√ьш\/} тхЁ°шэрьш (ўхЁхч ЁхсЁю $\mathbf{e}$) ш сєфхь
яшёрЄ№ $x\sim y$.
┬ё■фє эшцх
уЁрЇ $G$ яЁхфяюырурхЄё  ыюъры№эю ъюэхўэ√ь, Є.х.: т ы■сющ
юуЁрэшўхээющ юсырёЄш т $\rd$ ёюфхЁцшЄё  ъюэхўэюх ўшёыю тхЁ°шэ ш
ЁхсхЁ уЁрЇр $G$. ╩Ёюьх Єюую, уЁрЇ $G$ яЁхфяюырурхЄё 
$\zd$-яхЁшюфшўхёъшь:
\begin{equation}\label{KorSlo.1.1}
G+n=G,\ \ n\in\zd.
\end{equation}
═р уЁрЇх $G$ юяЁхфхышь фшёъЁхЄэ√щ юяхЁрЄюЁ ╦ряырёр $\Delta$,
фхщёЄтє■∙шщ т $\ell_{2}(X)$,
\begin{equation*}
\Delta u(x)=\sum_{y\sim x}(u(x)-u(y)),\ \ x\in X,\ \ u\in
\ell_{2}(X).
\end{equation*}
╧юфўхЁъэхь, ўЄю т яюёыхфэхщ ЇюЁьєых ёєььшЁютрэшх яЁюшёїюфшЄ яю
тёхь ЁхсЁрь шэЎшфхэЄэ√ь тхЁ°шэх $x$.

╬яхЁрЄюЁ $\Delta$ --- юуЁрэшўхээ√щ эхюЄЁшЎрЄхы№э√щ ёрьюёюяЁ цхээ√щ
юяхЁрЄюЁ т $\ell_{2}(X)$ (ёь., эряЁшьхЁ, \cite{KorSlo.KS2014}).
╬яЁхфхышь яхЁшюфшўхёъшщ юяхЁрЄюЁ ╪ЁхфшэухЁр $H=\Delta+Q$, уфх $Q$
хёЄ№ тх∙хёЄтхээ√щ яхЁшюфшўхёъшщ яюЄхэЎшры эр $X$, єфютыхЄтюЁ ■∙шщ
єёыютш■
\begin{equation}\label{KorSlo.1.2}
Q=\bar{Q}\in \ell_{\infty}(X);\ \ Q(x+n)=Q(x),\ \ (x,n)\in
X\times\zd.
\end{equation}

 \pr {\bf ╤яхъЄЁры№эюх Ёрчыюцхэшх эхтючьє∙хээюую
юяхЁрЄюЁр $H$.} ╬сючэрўшь ўхЁхч $\{x_{j}\}_{j=1}^{\nu}$
яхЁхёхўхэшх ьэюцхёЄтр тхЁ°шэ уЁрЇр $G$ ш  ўхщъш
$\Omega:=[0,1)^{d}$ Ёх°хЄъш $\zd$. ┬ ёшыє $\zd$-яхЁшюфшўэюёЄш
уЁрЇр $G$, фы  ърцфющ тхЁ°шэ√ $x\in X$ эрщфєЄё  хфшэёЄтхээ√х
$j\in\{1,\dots,\nu\}$ ш $n\in\zd$ Єръшх, ўЄю $x=x_{j}+n$. ┬тхфхь
єэшЄрЁэ√щ юяхЁрЄюЁ $\mathbb{U}:\ell_{2}(X)\to
\ell_{2}(\zd,\mathbb{C}^{\nu})$
\begin{equation}\label{KorSlo.1.3}
(\mathbb{U}u)_{j}(n):=u(x_{j}+n),\ \ j=1,\dots,\nu,\ \ n\in\zd,\ \
u\in \ell_{2}(X).
\end{equation}
╟рфрфшь т $\ell_{2}(\zd,\mathbb{C}^{\nu})$ юяхЁрЄюЁ
$\mathbb{H}:=\mathbb{U}H\mathbb{U}^{*}$ єэшЄрЁэю ¤ътштрыхэЄэ√щ
эхтючьє∙хээюьє юяхЁрЄюЁє $H$. ╬яЁхфхышь фшёъЁхЄэюх яЁхюсЁрчютрэшх
╘єЁ№х $$\Phi_{}:\ell_{2}(\zd,\mathbb{C}^{\nu})\to
L_{2}(\T^{d},\mathbb{C}^{\nu}),\ \ \T^{d}=[-\pi,\pi]^{d};$$
шьхээю, фы  тё ъющ $u\in \ell_{2}(\zd,\mathbb{C}^{\nu})$ яюыюцшь
\begin{equation}\label{KorSlo.1.4}
\left(\Phi_{}u\right)(x):=(2\pi)^{-d/2}\sum_{n\in\zd}e^{-inx}u(n),\
\ x\in\T^{d}.
\end{equation}
╙эшЄрЁэ√щ юяхЁрЄюЁ $\Phi_{}$ фшруюэрышчєхЄ (ёь., эряЁшьхЁ,
\cite{KorSlo.KS2014}) юяхЁрЄюЁ $\mathbb{H}$
\begin{equation}\label{KorSlo.1.5}
\mathbb{H}=\Phi_{}^{*}[h(k)]\Phi_{}.
\end{equation}
╟фхё№ $h(k)$, $k\in\T^{d}$, эхъюЄюЁр   ёрьюёюяЁ цхээр 
$\nu\times\nu$-ьрЄЁшЎр-ЇєэъЎш   рэрышЄшўхёър  яю $k$ (яюфЁюсэюх
юяшёрэшх $h(k)$, фрэю т \cite{KorSlo.KS2014}). ╬сючэрўшь ўхЁхч
$E_{1}(k)\le E_{2}(k)\le\dots\le E_{\nu}(k)$ ёюсёЄтхээ√х чэрўхэш 
ьрЄЁшЎ√ $h(k)$, чрэєьхЁютрээ√х т яюЁ фъх тючЁрёЄрэш , ё єўхЄюь
ъЁрЄэюёЄш. ╘єэъЎшш $E_{s}(k)$, $k\in\T^{d}$, $s=1,\dots,\nu$ сєфхь
эрч√трЄ№ чюээ√ьш ЇєэъЎш ьш. ╩рцфр  чюээр  ЇєэъЎш  $E_{s}(k)$ ---
эхяЁхЁ√тэр , ъєёюўэю тх∙хёЄтхээю-рэрышЄшўхёър  яю $k\in\T^{d}$.
╤яхъЄЁ $\sigma(H)=\sigma(\mathbb{H})$  ты хЄё  юс·хфшэхэшхь чюэ
$\sigma_s=E_{s}(\T^{d})$, ърцфр  шч ъюЄюЁ√ї хёЄ№ юсЁрч юЄюсЁрцхэш 
$k\mapsto E_{s}(k)$,
\begin{equation*}
\sigma(H)=\sigma(\mathbb{H})=\bigcup_{s=1}^{\nu}E_{s}(\T^{d}).
\end{equation*}
╟юэ√ ьюуєЄ яхЁхъЁ√трЄ№ё  ш ьюуєЄ т√ЁюцфрЄ№ё  т Єюўъє (яюёыхфэ  
юърч√трхЄё  ёюсёЄтхээ√ь чэрўхэшхь юяхЁрЄюЁр $H$ схёъюэхўэющ
ъЁрЄэюёЄш).  ╬ЄьхЄшь, ўЄю яхЁтр  чюэр эх т√ЁюцфрхЄё 
\cite{KorSlo.KS2014}. ╤яхъЄЁ $\sigma(H)$ ёюфхЁцшЄ фтх
яюыєсхёъюэхўэ√х ыръєэ√
\begin{equation}
\begin{aligned}
\label{KorSlo.1.6}
\sigma(H)\cap(-\infty,\Lambda_{\min})=\varnothing,\  \ \ \ \ \
\sigma(H)\cap(\Lambda_{\max},+\infty)=\varnothing,\\
\Lambda_{\min}:=\min\{E_{1}(k),k\in\T^{d}\},\ \
\Lambda_{\max}:=\max\{E_{\nu}(k),k\in\T^{d}\}.
\end{aligned}
\end{equation}
╩Ёюьх Єюую, ёяхъЄЁ юяхЁрЄюЁр $H$ ьюцхЄ ёюфхЁцрЄ№
тэєЄЁхээ■■ ыръєэє $(\Lambda_{+},\Lambda_{-})$, єфютыхЄтюЁ ■∙є■ яЁш
эхъюЄюЁюь $N\in\{2,\dots,\nu\}$ єёыютш■
\begin{equation}\label{KorSlo.1.7}
\Lambda_{+}=\max\{E_{N-1}(k),k\in\T^{d}\}<
\Lambda_{-}=\min\{E_{N}(k),k\in\T^{d}\}.
\end{equation}
\pr{\bf ╩ырёё√ ЇєэъЎшщ.}  ╧Ёхцфх, ўхь уютюЁшЄ№ ю тючьє∙р■∙хь
юяхЁрЄюЁх, ттхфхь эхюсїюфшь√х ъырёё√ ЇєэъЎшщ. ╧єёЄ№ $(\Z,d\mu)$
--- ёхярЁрсхы№эюх шчьхЁшьюх яЁюёЄЁрэёЄтю ё $\sigma$-ъюэхўэющ
ьхЁющ. ═рЁ фє ёю ёЄрэфрЁЄэ√ьш ъырёёрьш $L_{p}(\Z,d\mu)$ юяЁхфхышь
{\it ёырс√х\/} $L_{p}$-ъырёё√ $L_{p,\infty}(\Z,d\mu)$,
$p\in(0,+\infty)$, (ёь., эряЁшьхЁ,
\cite{KorSlo.BirmanSolomyakSnumbers77}, \cite{KorSlo.BKS91}).
╚ьхээю, фы  $\mu$-шчьхЁшьющ ЇєэъЎшш $f:\Z\to\mathbb{C}$ яюыюцшь
$$
\mu_{f}(s):=\mu(\{z\in\Z:|f(z)|>s\}),\ \ s>0.$$ ╩ырёё
$L_{p,\infty}$ т√фхы хЄё  ЄЁхсютрэшхь ъюэхўэюёЄш ЇєэъЎшюэрыр
\begin{equation}\label{KorSlo.1.8}
\|f\|_{L_{p,\infty}}:=\sup_{s>0}s\mu_{f}^{1/p}(s).
\end{equation}
╧ЁюёЄЁрэёЄтю $L_{p,\infty}$ яюыэю юЄэюёшЄхы№эю ътрчшэюЁь√
(\ref{KorSlo.1.8}), 
ш ёюфхЁцшЄ
ёхярЁрсхы№эюх яюфяЁюёЄЁрэёЄтю
\begin{equation*}
L_{p,\,\infty}^{0}:=\{f\in L_{p,\infty}:\mu_{f}(s)=o(s^{-p}),\ \
s\to+0,\ \ s\to+\infty\}.
\end{equation*}
╧ЁюёЄЁрэёЄтю $L_{p,\infty}(\Z,d\mu)$ ёхярЁрсхы№эю, хёыш ш Єюы№ъю
хёыш $\Z$  ты хЄё  юс·хфшэхэшхь ъюэхўэюую ўшёыр рЄюьют; т ¤Єюь
ёыєўрх $L_{p,\infty}=L_{p,\infty}^{0}$ --- ъюэхўэюьхЁэюх. ┼ёыш
$\Z$ --- шчьхЁшьюх яюфьэюцхёЄтю т $\rd$, ш $d\mu$ --- ьхЁр ╦хсхур,
Єю $L_{p,\infty}(\Z,d\mu)=:L_{p,\infty}(\Z)$.

╧Ёштхфхь эхёъюы№ъю яЁшьхЁют ЇєэъЎшщ, яЁшэрфыхцр∙шї ёырс√ь
$L_{p}$-ъырёёрь (ёь., эряЁшьхЁ, \cite{KorSlo.BirmanSolomyak92}, р
Єръцх \cite{KorSlo.BirmanSolomyak79}).
\begin{Ex}
{\rm 1)} ╧єёЄ№ ЇєэъЎш 
$V(x)=\vartheta(\frac{x}{|x|})\cdot|x|^{-d/p}$, $x\in\rd$, уфх
$\vartheta\in L_{p}(\mathbb{S}^{d-1})$ яЁш эхъюЄюЁюь
$p\in(0,+\infty)$. ╥юуфр ёяЁртхфышт√ тъы■ўхэшх $V\in
L_{p,\infty}(\rd)$ ш ЁртхэёЄтю
$\|V\|_{L_{p,\infty}}=d\,{}^{-\frac{1}{p}}\|\vartheta\|_{L_{p}}$.

\noindent{\rm 2)} ┼ёыш ЇєэъЎш  $V\in L_{\infty}(\rd)$
єфютыхЄтюЁ хЄ єёыютш■ $V(x)=O(|x|^{-\frac{d}{p}})$, яЁш
$|x|\to\infty$, Єю ёяЁртхфыштю тъы■ўхэшх $V\in L_{p,\infty}(\rd)$.

\noindent{\rm 3)} ┼ёыш ЇєэъЎш  $V\in L_{\infty}(\rd)$ єфютыхЄтюЁ хЄ
єёыютш■ $V(x)=o(|x|^{-\frac{d}{p}})$,  яЁш $|x|\to\infty$, Єю
ёяЁртхфыштю тъы■ўхэшх $V\in L_{p,\infty}^{0}(\rd)$.
\end{Ex}
╧Ёхфяюыюцшь, ЄхяхЁ№, ўЄю эр $\Z$ чрфрэр $\mu$-шчьхЁшьр 
ьрЄЁшўэючэрўэр  ЇєэъЎш  $\W$. ┬ ¤Єюь ёыєўрх яюыюцшь яю юяЁхфхыхэш■
$$
\W\in L_{p,\infty}(\Z,d\mu)\Longleftrightarrow\forall i,j\ \
\W_{i,j}\in L_{p,\infty}(\Z,d\mu);\ \
\|\W\|_{L_{p,\infty}}=\max_{i,j}\|\W_{i,j}\|_{L_{p,\infty}}.
$$
┼ёыш $\Z$ --- ёўхЄэюх ьэюцхёЄтю, ь√ сєфхь яюы№чютрЄ№ё 
юсючэрўхэш ьш
$$L_{p}(\Z,d\mu)=:\ell_{p}(\Z,d\mu),\ \ L_{p,\infty}(\Z,d\mu)=:\ell_{p,\infty}(\Z,d\mu);$$
хёыш фюяюыэшЄхы№эю $d\mu$ --- ёўшЄр■∙р  ьхЁр, Єю
$$\ell_{p}(\Z,d\mu)=:\ell_{p}(\Z),\ \ \ell_{p,\infty}(\Z,d\mu)=:\ell_{p,\infty}(\Z).$$
╧ютЄюЁ   Ёрёёєцфхэш  шч \cite{KorSlo.BirmanSolomyak92}, яюыєўшь
фы  $m_{1}\times m_{2}$-ьрЄЁшЎ√-ЇєэъЎшш $\W(n)$, $n\in\zd$,
ёыхфє■∙хх єЄтхЁцфхэшх.
\begin{PPP}\label{KorSlo.Note1.1}
{\rm 1)} ┼ёыш $\W(n)=O(|n|^{-d/p})$, Єю $\W\in
\ell_{p,\infty}(\zd)$.

{\rm 2)} ┼ёыш $\W(n)=o(|n|^{-d/p})$, Єю $\W\in
\ell_{p,\infty}^{\,0}(\zd)$.

{\rm 3)} ┼ёыш $\W(0)=0$ ш $\W(n)=\omega(\frac{n}{|n|})\cdot
|n|^{-d/p}$, $n\not=0$, уфх  $\omega\in B(\mathbb{S}^{d-1})$, Єю
$\W\in \ell_{p,\infty}(\zd)$, $\|\W\|_{\ell_{p,\infty}}\le
C\|\omega\|_{B}$, уфх $C=C(m_{1},m_{2},p,d)$.
%
\end{PPP}
\pr{\bf ┬ючьє∙хэшх.} ╠√ ЁрёёьюЄЁшь тючьє∙хээ√щ юяхЁрЄюЁ ╪ЁхфшэухЁр
$H_{\pm}(t)=H\pm tV$, $t>0$, уфх
яЁхфяюырурхЄё , ўЄю  яюЄхэЎшры $V$ чэръююяЁхфхыхэ, юуЁрэшўхэ ш
єс√трхЄ эр схёъюэхўэюёЄш
\begin{equation}\label{KorSlo.1.9}
0\le V\in\ell_{\infty}(X),\ \ V(x)\to 0,\ \ |x|\to+\infty.
\end{equation}
┴юыхх Єюўэю, ь√ ЁрёёьюЄЁшь тючьє∙хэш  єфютыхЄтюЁ ■∙шх юфэюьє шч
ЄЁхї єёыютшщ:
\begin{equation}\label{KorSlo.1.10}
0\le V\in\ell_{p}(X);
\end{equation}
\begin{equation}\label{KorSlo.1.11}
0\le V\in\ell_{p,\infty}(X);
\end{equation}
\begin{equation}\label{KorSlo.1.12}
0\le V\in \ell_{\infty}(X),\ \
V(x)=|x|^{-d/p}\left(\vartheta\left(\frac{x}{|x|}\right)+o(1)\right),\
\ x\in X,\ \ |x|\to\infty,\ \ \vartheta\in C(\mathbb{S}^{d-1}).
\end{equation}
╬ЄьхЄшь, ўЄю яЁш єёыютшш (\ref{KorSlo.1.9}) (р чэрўшЄ ш яЁш ы■сюь шч
єёыютшщ (\ref{KorSlo.1.10}) -- (\ref{KorSlo.1.12})) юяхЁрЄюЁ
єьэюцхэш  эр $V$ ъюьяръЄхэ. ═рё шэЄхЁхёєхЄ фшёъЁхЄэ√щ ёяхъЄЁ
юяхЁрЄюЁют $H_{\pm}(t)=H\pm tV$, $t>0$. ╧Ёхфяюыюцшь, ўЄю
$(\Lambda_{+},\Lambda_{-})$ --- ыръєэр т ёяхъЄЁх юяхЁрЄюЁр $H$
(тючьюцэю яюыєсхёъюэхўэр ). ╧юёъюы№ъє тючьє∙хэшх $V$
--- ъюьяръЄэ√щ юяхЁрЄюЁ т $\ell_{2}(X)$, ёяхъЄЁ√ тючьє∙хээ√ї
юяхЁрЄюЁют $H_{\pm}(t)$, $t>0$, т ыръєэх
$(\Lambda_{+},\Lambda_{-})$ фшёъЁхЄэ√. ╬ёэютэ√ь юс·хъЄюь
шёёыхфютрэш   ты ■Єё  ёўшЄр■∙шх ЇєэъЎшш
\begin{equation}\label{KorSlo.1.13}
N_{\pm}(\lambda,\tau):=\sum_{t\in(0,\tau)}\mathrm{dimKer}(H_{\pm}(t)-\lambda
I),\ \ \tau>0,\ \ \lambda\in(\Lambda_{+},\Lambda_{-}).
\end{equation}
╧Ёртр  ўрёЄ№ т (\ref{KorSlo.1.13}) ёўшЄрхЄё  схёъюэхўэющ, хёыш
эрщфхЄё  схёъюэхўэюх ўшёыю Єюўхъ $t\in(0,\tau)$ т ъюЄюЁ√ї
$\mathrm{dimKer}(H_{\pm}(t)-\lambda I)>0$. ╚ч трЁшрЎшюээ√ї
ёююсЁрцхэшщ т√ЄхърхЄ ёыхфє■∙хх єЄтхЁцфхэшх (яЁшэЎшя
┴шЁьрэр-╪тшэухЁр, ёь., эряЁшьхЁ, \cite{KorSlo.Birman91}).
\begin{PPP}\label{KorSloPropositionBirmanShwinger}
╧єёЄ№ $H$ --- ёрьюёюяЁ цхээ√щ юяхЁрЄюЁ т $\ell_{2}(X)$, $V$ ---
эхюЄЁшЎрЄхы№э√щ ъюьяръЄэ√щ юяхЁрЄюЁ т $\ell_{2}(X)$,
$\sigma(H)\cap(\Lambda_{+},\Lambda_{-})=\varnothing$. ╥юуфр
ёяЁртхфыштю ЁртхэёЄтю
\begin{equation}\label{KorSlo.1.14}
N_{\pm}(\lambda,\tau)=n_{\pm}(\tau^{-1},V^{1/2}(\lambda
I-H)^{-1}V^{1/2}),\ \ \tau>0,\ \
\lambda\in(\Lambda_{+},\Lambda_{-}).
\end{equation}
\end{PPP}
╚ч яЁхфыюцхэш  \ref{KorSloPropositionBirmanShwinger} т√ЄхърхЄ, ўЄю
тхышўшэ√ $N_{\pm}(\lambda,\tau)$, $\tau>0$,
$\lambda\in(\Lambda_{+},\Lambda_{-})$, ъюэхўэ√ ш ьюэюЄюээ√ яю
$\lambda$. ╤ыхфютрЄхы№эю, $N_{\pm}(\lambda,\tau)$ ьюцэю юяЁхфхышЄ№ ш
яЁш $\lambda=\Lambda_{\pm}$, ъръ яЁхфхы (тючьюцэю схёъюэхўэ√щ)
\begin{equation*}
N_{\pm}(\Lambda_{\pm},\tau):=\lim_{\lambda\to\Lambda_{\pm}\pm
0}N_{\pm}(\lambda,\tau),\ \ \tau>0.
\end{equation*}
\par\noindent{\pr\bf ╬ёэютэющ Ёхчєы№ЄрЄ.}

{\it ═р яЁюЄ цхэшш ¤Єюую яєэъЄр ь√ яЁхфяюырурхь т√яюыэхээ√ьш
єёыютш  яхЁшюфшўэюёЄш {\rm(\ref{KorSlo.1.1})} ш
{\rm(\ref{KorSlo.1.2})}; ъЁюьх Єюую яЁхфяюырурхЄё , ўЄю
$(\Lambda_{+},\Lambda_{-})$
--- ыръєэр т ёяхъЄЁх юяхЁрЄюЁр $H$ (тючьюцэю яюыєсхёъюэхўэр ); эхюЄЁшЎрЄхы№э√щ яюЄхэЎшры $V$ єфютыхЄтюЁ хЄ юфэюьє шч єёыютшщ {\rm(\ref{KorSlo.1.10})--(\ref{KorSlo.1.12})}\/}.

═шцх ь√ яюърцхь, ўЄю яЁш єёыютш ї (\ref{KorSlo.1.10}) ышсю
(\ref{KorSlo.1.11}) ш эхъюЄюЁ√ї яЁхфяюыюцхэш ї ю яютхфхэшш чюээ√ї
ЇєэъЎшщ ёяхъЄЁ тючьє∙хээюую юяхЁрЄюЁр т ыръєэх
$(\Lambda_{+},\Lambda_{-})$ ъюэхўхэ; сєфхЄ фрэр юЎхэър тхышўшэ
$N_{\pm}(\Lambda_{\pm},\tau)$. ╠√ яюърцхь, ўЄю яЁш єёыютшш
(\ref{KorSlo.1.12}) ёяЁртхфышт√ рёшьяЄюЄшъш (\ref{KorSlo.0.2}), ш
рёшьяЄюЄшўхёъшх ъю¤ЇЇшЎшхэЄ√
$\Gamma_{p}^{\pm}(\lambda,H,V)=\Gamma_{p}^{\pm}(\lambda)$
юяЁхфхы ■Єё  ЁртхэёЄтрьш
\begin{equation}\label{KorSlo.1.15}
\Gamma_{p}^{\pm}(\lambda):=
\frac{1}{d(2\pi)^{d}}\sum_{s=1}^{\nu}\int_{\T^{d}}\left(\lambda-E_{s}(k)\right)_{\pm}^{-p}dk
\int_{\mathbb{S}^{d-1}}\vartheta^{p}(\theta)dS(\theta),\ \
\lambda\in[\Lambda_{+},\Lambda_{-}].
\end{equation}
╧Ёш тёхї $\lambda\in(\Lambda_{+},\Lambda_{-})$ тхышўшэ√
$\Gamma_{p}^{\pm}(\lambda)$ ъюэхўэ√. ╬ЄьхЄшь, ўЄю тхышўшэ√
$\Gamma_{p}^{\pm}(\Lambda_{\pm})$ ъюэхўэ√, хёыш т√яюыэхэ√
(ёююЄтхЄёЄтхээю) єёыютш 
\begin{equation}
\label{KorSlo.1.16} (\Lambda_{\pm}-E_{s}(\cdot))_{\pm}^{-1}\in
L_{\varkappa}(\T^{d}),\ \ s=1,\dots,\nu,\ \ \text{уфх}\ \ \left\{
\begin{array}{l}
\varkappa=p,\ \ p>1,\\
\varkappa=1,\ \ p\in(0,1),\\
\varkappa>1,\ \ p=1.
\end{array}
\right.
\end{equation}
╧Ёш $p>1$ ъЁюьх єёыютшщ (\ref{KorSlo.1.16}) эрь яюэрфюс Єё  ш
сюыхх ёырс√х єёыютш 
\begin{equation}\label{KorSlo.1.17}
(\Lambda_{\pm}-E_{s}(\cdot))_{\pm}^{-1}\in L_{p,\infty}(\T^{d}),\
\ s=1,\dots,\nu.
\end{equation}
╤ыхфє■∙шх фтх ЄхюЁхь√ фр■Є юЎхэъш ёўшЄр■∙хщ ЇєэъЎшш
$N_{\pm}(\Lambda_{\pm},\tau)$, $\tau>0$.
\begin{Th}\label{KorSloTheoremRLC}
╧єёЄ№
яЁш эхъюЄюЁюь
$p\in(1,+\infty)$ т√яюыэхэ√ єёыютш  {\rm(\ref{KorSlo.1.10})} ш
{\rm(\ref{KorSlo.1.17}$\pm$)}. ╥юуфр ёяхъЄЁ юяхЁрЄюЁр
$H_{\pm}(t)$, $t>0$, т ыръєэх $(\Lambda_{+},\Lambda_{-})$ ъюэхўхэ
ш ёяЁртхфышт√ юЎхэъш
\begin{equation}\label{KorSlo.1.18}
N_{\pm}(\Lambda_{\pm},\tau)\le
C(p,d,\nu)\,\tau^{p}\sum_{s=1}^{\nu}\|(\Lambda_{\pm}-E_{s}(\cdot))_{\pm}^{-1}\|_{L_{p,\infty}}^{p}\|V\|_{\ell_{p}}^{p},\
\ \tau>0;
\end{equation}
\begin{equation}\label{KorSlo.1.19}
N_{\pm}(\Lambda_{\pm},\tau)=o(\tau^{p}),\ \ \tau\to+\infty.
\end{equation}
\end{Th}
╟фхё№ ш фрыхх ЇюЁьєы√ т ЄхъёЄх ўрёЄю ёэрсцхэ√ шэфхъёрьш <<$+$>> ш
<<$-$>>. ┬ ¤Єюь ёыєўрх (хёыш юёюсю эх юуютюЁхэю фЁєуюх) шї ёыхфєхЄ
ўшЄрЄ№ эхчртшёшью фы  ърцфюую шч фтєї шэфхъёют. ┼ёыш ь√ їюЄшь
ёяхЎшры№эю т√фхышЄ№ ёыєўрщ юфэюую чэрър, ь√ єърч√трхь хую яЁш
ёёы√ыъх эр ёююЄтхЄёЄтє■∙є■ ЇюЁьєыє, эряЁшьхЁ єёыютшх
(\ref{KorSlo.1.16}$+$) ш Є.я.
\begin{Th}\label{KorSloTheoremWRLC}
╧єёЄ№
яЁш эхъюЄюЁюь $p\in(0,+\infty)$
т√яюыэхэ√ єёыютш  {\rm(\ref{KorSlo.1.11})} ш
{\rm(\ref{KorSlo.1.16}$\pm$)}. ╥юуфр ёяхъЄЁ юяхЁрЄюЁр
$H_{\pm}(t)$, $t>0$, т ыръєэх $(\Lambda_{+},\Lambda_{-})$ ъюэхўхэ
ш ёяЁртхфыштр юЎхэър
\begin{equation}\label{KorSlo.1.20}
N_{\pm}(\Lambda_{\pm},\tau)\le
C(p,\varkappa,d,\nu)\,\tau^{p}\sum_{s=1}^{\nu}\|(\Lambda_{\pm}-E_{s}(\cdot))_{\pm}^{-1}\|_{L_{\varkappa}}^{p}
\|V\|_{\ell_{p,\infty}}^{p},\ \ \tau>0.
\end{equation}
┼ёыш фюяюыэшЄхы№эю $V\in\ell_{p,\infty}^{\,0}(X)$, Єю
$N_{\pm}(\Lambda_{\pm},\tau)=o(\tau^{p})$, $\tau\to+\infty$.
\end{Th}

╬ёэютэ√ь Ёхчєы№ЄрЄюь ЁрсюЄ√  ты хЄё  ёыхфє■∙р  ЄхюЁхьр.
\begin{Th}\label{KorSloTheoremAsymptotic}
╧єёЄ№
яЁш эхъюЄюЁюь
$p\in(0,+\infty)$ т√яюыэхэю єёыютшх {\rm(\ref{KorSlo.1.12})}.
╥юуфр яЁш тёхї $\lambda\in(\Lambda_{+},\Lambda_{-})$ шьх■Є ьхёЄю
рёшьяЄюЄшъш
\begin{equation}
\label{KorSlo.1.21}
N_{\pm}(\lambda,\tau)=\tau^{p}(\Gamma_{p}^{\pm}(\lambda)+o(1)),\ \
\tau\to+\infty.
\end{equation}
┼ёыш, ъЁюьх Єюую, т√яюыэхэю єёыютшх {\rm(\ref{KorSlo.1.16}$\pm$)},
Єю рёшьяЄюЄшър {\rm(\ref{KorSlo.1.21}$\pm$)} тхЁэр ш яЁш
$\lambda=\Lambda_{\pm}$.
\end{Th}
╧Ёш эхэєыхтюь тючьє∙хэшш ъю¤ЇЇшЎшхэЄ $\Gamma^{+}_{p}(\lambda)$
яюыюцшЄхыхэ, хёыш $\lambda$ ыхцшЄ тю тэєЄЁхээхщ ыръєэх шыш т
яЁртющ яюыєсхёъюэхўэющ ыръєэх; рэрыюушўэю ъю¤ЇЇшЎшхэЄ
$\Gamma^{-}_{p}(\lambda)$ яюыюцшЄхыхэ, хёыш $\lambda$ ыхцшЄ тю
тэєЄЁхээхщ ыръєэх шыш т ыхтющ яюыєсхёъюэхўэющ ыръєэх. ╨рчєьххЄё 
$N_{+}(\lambda,\tau)=\Gamma_{p}^{+}(\lambda)=0$ яЁш
$\lambda\le\Lambda_{\min}$, ш
$N_{-}(\lambda,\tau)=\Gamma_{p}^{-}(\lambda)=0$ яЁш
$\lambda\ge\Lambda_{\max}$.

╧ютхфхэшх чюэ√ї ЇєэъЎшщ  $E_s(k)$ эр ъЁр■ ыръєэ√ шчєўрыюё№ т
\cite{KorSlo.KS2016}, \cite{KorSlo.KS2018}.  ═рь яюЄЁхсєхЄё 
ёыхфє■∙хх юяЁхфхыхэшх Ёхуєы ЁэюёЄш ъЁр  ыръєэ√.

{\it 
╦хт√щ ъЁрщ ыръєэ√ $(\Lambda_{+},\Lambda_{-})$ эрчютхь Ёхуєы Ёэ√ь,
хёыш Єюўър $\Lambda_{+}$ ёюфхЁцшЄё  т юсЁрчх Єюы№ъю юфэющ чюээющ
ЇєэъЎшш $E_{N-1}(k)$, $k\in\T^{d}$, ш ЁртхэёЄтю
$E_{N-1}(k)=\Lambda_{+}$ фюёЄшурхЄё  т ъюэхўэюь ўшёых Єюўхъ
$k_{j}\in \T^{d}$, $j=1,...,M$, ърцфр  шч ъюЄюЁ√ї хёЄ№
эхт√Ёюцфхээр  Єюўър ьръёшьєьр фы  $E_{N-1}(\cdot)$, Є.х.
\begin{equation*}
\Lambda_{+} -E_{N-1}(k)=q_{j}(k-k_{j})+O(|k-k_{j}|^{3}),\ \ k\to
k_{j},\ \ j=1,...,M,
\end{equation*}
уфх $q_{j}$ --- яюыюцшЄхы№эю юяЁхфхыхээр  ътрфЁрЄшўэр  ЇюЁьр.
└эрыюушўэ√ь юсЁрчюь юяЁхфхы хЄё  Ёхуєы ЁэюёЄ№ яЁртюую ъЁр  ыръєэ√.
\/}

┼ёыш ъЁрщ ыръєэ√ $\Lambda_{\pm}$ Ёхуєы Ёхэ, Єю
\begin{gather*}
(\Lambda_{\pm}-E_{s}(\cdot))_{\pm}^{-1}\in L_{\varkappa}(\T^{d}),\
\
s=1,...,\nu,\ \ \forall  \varkappa< \frac d2;\\
(\Lambda_{\pm}-E_{s}(\cdot))_{\pm}^{-1}\in
L_{\varkappa,\infty}(\T^{d}),\ \ s=1,...,\nu,\ \ \forall
\varkappa\le \frac d2.
\end{gather*}
╥ръшь юсЁрчюь, шч ЄхюЁхь \ref{KorSloTheoremRLC} ш
\ref{KorSloTheoremWRLC} ёыхфєхЄ єЄтхЁцфхэшх.
\begin{Con}\label{KorSlo.Corollary1.6}
{\rm 1)} ╧єёЄ№ $d\ge 3$ ш яЁш эхъюЄюЁюь $p\in(0,\frac{d}{2})$
т√яюыэхэю
єёыютшх {\rm(\ref{KorSlo.1.11})}, 
ш ъЁрщ ыръєэ√ $\Lambda_{\pm}$ Ёхуєы Ёхэ. ╥юуфр ёяхъЄЁ юяхЁрЄюЁр
$H_{\pm}(t)$, $t>0$, т ыръєэх $(\Lambda_{+},\Lambda_{-})$ ъюэхўхэ
ш ёяЁртхфыштр юЎхэър
%
\begin{equation*}
N_{\pm}(\Lambda_{\pm},\tau)\le
C_{\pm}(p,d,H,G)\|V\|_{\ell_{p,\infty}}^{p}\tau^{p},\ \ \tau>0.
\end{equation*}
┼ёыш фюяюыэшЄхы№эю $V\in\ell_{p,\infty}^{\,0}(X)$, Єю
$N_{\pm}(\Lambda_{\pm},\tau)=o(\tau^{p})$, $\tau\to+\infty$.
\par\noindent{\rm 2)} ╧єёЄ№ т√яюыэхэю єёыютшх 
$V\in\ell_{d/2}(X)$, $d\ge 3$,
ш ъЁрщ ыръєэ√ $\Lambda_{\pm}$
Ёхуєы Ёхэ. ╥юуфр ёяхъЄЁ юяхЁрЄюЁр $H_{\pm}(t)$, $t>0$, т ыръєэх
$(\Lambda_{+},\Lambda_{-})$ ъюэхўхэ ш ёяЁртхфышт√ юЎхэъш
\begin{gather*}
N_{\pm}(\Lambda_{\pm},\tau)\le
C_{\pm}(p,d,H,G)\|V\|_{\ell_{d/2}}^{d/2}\tau^{d/2},\ \ \tau>0;\\
N_{\pm}(\Lambda_{\pm},\tau)=o(\tau^{d/2}),\ \ \tau\to+\infty.
\end{gather*}
\end{Con}
└эрыюушўэю, шч ЄхюЁхь√ \ref{KorSloTheoremAsymptotic} т√ЄхърхЄ
ёыхфє■∙хх єЄтхЁцфхэшх.
\begin{Con}\label{KorSlo.Corollary1.7}
╧єёЄ№ яЁш эхъюЄюЁ√ї $d\ge 3$, $p\in(0,\frac d2)$, т√яюыэхэю
єёыютшх {\rm(\ref{KorSlo.1.12})}, 
ш ыхт√щ {\rm(}шыш яЁрт√щ{\rm)} ъЁрщ
ыръєэ√ Ёхуєы Ёхэ.
%
╥юуфр рёшьяЄюЄшър {\rm(\ref{KorSlo.1.21}$+$)}
{\rm(}шыш {\rm(\ref{KorSlo.1.21}$-$))} ёяЁртхфыштр яЁш
$\lambda=\Lambda_{+}$ {\rm(}яЁш $\lambda=\Lambda_{-}${\rm)}.
\end{Con}
╚ч Ёхчєы№ЄрЄют ЁрсюЄ ╩юЁюЄ хтр, ╤рсєЁютющ (\cite{KorSlo.KS2016},
ЄхюЁхьр 1.2; \cite{KorSlo.KS2018}, ЄхюЁхьр 2.1) ёыхфєхЄ, ўЄю яЁш
єёыютш ї (\ref{KorSlo.1.1}), (\ref{KorSlo.1.2}) эшцэшщ ъЁрщ
ёяхъЄЁр $\Lambda_{\min}$ юяхЁрЄюЁр $H$ Ёхуєы Ёхэ. ╥ръшь юсЁрчюь,
шч ёыхфёЄтшщ \ref{KorSlo.Corollary1.6} ш \ref{KorSlo.Corollary1.7}
т√ЄхърхЄ ёыхфє■∙хх єЄтхЁцфхэшх.
\begin{Th}\label{KorSloTheoremLowerEdge}
{\rm 1)} ╧єёЄ№ яЁш эхъюЄюЁ√ї $d\ge 3$, $p\in(0,\frac{d}{2})$
т√яюыэхэю єёыютшх {\rm(\ref{KorSlo.1.11})}. ╥юуфр ёяхъЄЁ юяхЁрЄюЁр
$H_{-}(t)$, $t>0$, ыхтхх Єюўъш $\Lambda_{\min}$ ъюэхўхэ ш
ёяЁртхфыштр юЎхэър
\begin{equation}\label{KorSlo.1.22}
N_{-}(\Lambda_{\min},\tau)\le
C(p,d,H,G)\|V\|_{\ell_{p,\infty}}^{p}\tau^{p},\ \ \tau>0.
\end{equation}
┼ёыш фюяюыэшЄхы№эю $V\in\ell_{p,\infty}^{0}(X)$, Єю
$N_{-}(\Lambda_{\min},\tau)=o(\tau^{p})$, яЁш $\tau\to+\infty$.
\par\noindent
{\rm 2)} ┼ёыш $V\in\ell_{d/2}(X)$, $d\ge 3$, Єю ёяхъЄЁ юяхЁрЄюЁр
$H_{-}(t)$, $t>0$, ыхтхх Єюўъш $\Lambda_{\min}$ ъюэхўхэ ш
ёяЁртхфышт√ юЎхэъш
\begin{equation}\label{KorSlo.1.23}
N_{-}(\Lambda_{\min},\tau)\le
C(p,d,H,G)\|V\|_{\ell_{d/2}}^{d/2}\tau^{d/2},\ \ \tau>0;
\end{equation}
\begin{equation}\label{KorSlo.1.24}
 N_{-}(\Lambda_{\min},\tau)=o(\tau^{d/2}),\ \
\tau\to+\infty.
\end{equation}
\noindent {\rm 3)} ╧єёЄ№ яЁш эхъюЄюЁ√ї $d\ge 3$,
$p\in(0,\frac{d}{2})$ т√яюыэхэю єёыютшх {\rm(\ref{KorSlo.1.12})}.
╥юуфр ёяхъЄЁ юяхЁрЄюЁр $H_{-}(t)$, $t>0$, ыхтхх Єюўъш
$\Lambda_{\min}$ ъюэхўхэ ш рёшьяЄюЄшър {\rm(\ref{KorSlo.1.21}$-$)}
ёяЁртхфыштр яЁш тёхї $\lambda\le\Lambda_{\min}$.
\end{Th}

\noindent{\bf ╩юььхэЄрЁшш.} $\bullet$ ╥хюЁхь√
\ref{KorSloTheoremRLC}--\ref{KorSloTheoremAsymptotic} ш ёыхфёЄтш 
\ref{KorSlo.Corollary1.6}, \ref{KorSlo.Corollary1.7} ёяЁртхфышт√
фы  яЁюшчтюы№эюую юяхЁрЄюЁр
$H=\mathbb{U}^{*}\Phi^{*}[h(k)]\Phi\mathbb{U}$, уфх $h(k)$ ---
ёрьюёюяЁ цхээр  $\nu\times\nu$-ьрЄЁшЎр-ЇєэъЎш  эхяЁхЁ√тэр  яю
$k\in\T^{d}$.
%

\noindent $\bullet$ ╬Ўхэъш (\ref{KorSlo.1.22})-(\ref{KorSlo.1.24})
т ёыєўрх $G=\zd$, $H=\Delta$ с√ыш яюыєўхэ√ т ЁрсюЄрї ╦хтшэр,
╨ючхэсы■ьр ш ╤юыюь ър \cite{KorSlo.LS97},
\cite{KorSlo.RozenblumSolomyak2008},
\cite{KorSlo.RozenblumSolomyak2009}, р Єръцх т ЁрсюЄх ┴рїр
\cite{KorSlo.Bach}. ┬ ЁрсюЄх ╨ючхэсы■ьр, ╤юыюь ър
\cite{KorSlo.RozenblumSolomyak2010} юЎхэъш
(\ref{KorSlo.1.22})-(\ref{KorSlo.1.24}) с√ыш яюыєўхэ√ фы  юс∙хую
уЁрЇр $G$ ш тхёютюую юяхЁрЄюЁр ╦ряырёр $H$ эр уЁрЇх $G$. ├ЁрЇ $G$
эх яЁхфяюырурыё  яхЁшюфшўхёъшь; ЄЁхсютрыюё№, ўЄюс√  фЁю юяхЁрЄюЁр
$e^{-tH}$, $t>0$, єфютыхЄтюЁ ыю юяЁхфхыхээ√ь єёыютш ь яЁш сюы№°шї
$t$.

\noindent $\bullet$ ╧юЁ фюъ рёшьяЄюЄшъш (\ref{KorSlo.1.21})
юЄышўрхЄё  юЄ яюЁ фър рёшьяЄюЄшъш (\ref{KorSlo.0.5}).
└ёшьяЄюЄшўхёъшщ ъю¤ЇЇшЎшхэЄ чртшёшЄ юЄ Єюўъш эрсы■фхэш  $\lambda$.
<<┬хщыхтёър >> яЁшЁюфр рёшьяЄюЄшъш (\ref{KorSlo.1.21}) ёЄрэютшЄё 
 ёэр, хёыш яюьхэ Є№ Ёюы ьш ъююЁфшэрЄ√ ш ътрчшшьяєы№ё√.

\noindent $\bullet$ ╨рчєьххЄё , фы  яюыєсхёъюэхўэ√ї ыръєэ
(\ref{KorSlo.1.6}) ёєььрЁэр  ъЁрЄэюёЄ№ ёяхъЄЁр юяхЁрЄюЁр
$H_{-}(\tau)$ т ыръєэх $(-\infty,\Lambda_{\min})$ ёютярфрхЄ ё
тхышўшэющ $N_{-}(\Lambda_{\min},\tau)$; ёєььрЁэр  ъЁрЄэюёЄ№
ёяхъЄЁр юяхЁрЄюЁр $H_{+}(\tau)$ т ыръєэх
$(\Lambda_{\max},+\infty)$ ёютярфрхЄ ё тхышўшэющ
$N_{+}(\Lambda_{\max},\tau)$.

\noindent $\bullet$ ╧Ёш $d=2$ т ёыєўрх Ёхуєы Ёэюую ыхтюую (шыш
яЁртюую) ъЁр  ыръєэ√ ЄхюЁхь√ \ref{KorSloTheoremRLC},
\ref{KorSloTheoremWRLC}, р Єръцх ЄхюЁхьр
\ref{KorSloTheoremAsymptotic} яЁш $\lambda=\Lambda_{\pm}$ эшўхую эх
фр■Є, Є.ъ. єёыютш  {\rm(\ref{KorSlo.1.16})},
{\rm(\ref{KorSlo.1.17})} чртхфюью эх т√яюыэхэ√.

\noindent $\bullet$ ╧юЁ фюъ рёшьяЄюЄшъш (\ref{KorSlo.1.21}) тю {\it
тэєЄЁхээхщ\/} Єюўъх ёяхъЄЁры№эющ ыръєэ√ юяхЁрЄюЁр $H$ эх чртшёшЄ юЄ
чэрър тючьє∙хэш .

\noindent $\bullet$ ╬яхЁрЄюЁ√ $H_{+}(t)$, $H_{-}(t)$, $t>0$,
тїюф Є ёшььхЄЁшўэю т ЄхюЁхь√ \ref{KorSloTheoremRLC},
\ref{KorSloTheoremWRLC} ш \ref{KorSloTheoremAsymptotic}. ╬фэръю
єёЄЁющёЄтю ыхтюую ш яЁртюую ъЁрхт ёяхъЄЁры№эющ ыръєэ√ юяхЁрЄюЁр
$H$ ьюцхЄ с√Є№ Ёрчэ√ь, р яюЄюьє яютхфхэшх тхышўшэ√
$N_{+}(\Lambda_{+},\tau)$ ьюцхЄ, тююс∙х уютюЁ , юЄышўрЄ№ё  юЄ
яютхфхэш  тхышўшэ√ $N_{-}(\Lambda_{-},\tau)$.

\noindent $\bullet$ ┬ ёыєўрх {\it фтєфюы№эюую\/} уЁрЇр (уЁрЇр,
тхЁ°шэ√ ъюЄюЁюую Ёрёярфр■Єё  эр фтр эхяхЁхёхър■∙шїё  ъырёёр, яЁшўхь
ЁхсЁр ёюхфшэ ■Є ыш°№ тхЁ°шэ√ Ёрчэ√ї ъырёёют) ёяхъЄЁ юяхЁрЄюЁр
╦ряырёр $\Delta$ ёшььхЄЁшўхэ юЄэюёшЄхы№эю хую ёхЁхфшэ√, ш яЁрт√щ
ъЁрщ ёяхъЄЁр юяхЁрЄюЁр $H=\Delta$ Єръцх Ёхуєы Ёхэ (ёь.
\cite{KorSlo.M82}, \cite{KorSlo.KS2016}) . ┬ ¤Єюь ёыєўрх ёяЁртхфышт
рэрыюу ЄхюЁхь√ \ref{KorSloTheoremLowerEdge} фы  яЁртюую ъЁр  ыръєэ√.

\section{╧ЁхфтрЁшЄхы№э√х ётхфхэш }

╧єёЄ№ чрфрэ√ ьрЄЁшЎ√-ЇєэъЎшш $f(k)$, $\V(n)$, $g(k)$, $n\in\zd$,
$k\in\T^{d}:=[-\pi,\pi]^{d}$, ёюуырёютрээ√ї яюЁ фъют ё Єхь, ўЄюс√
яЁюшчтхфхэшх $\mathbb{G}(k,n):=f(k)\V(n)g(k)$ яЁхфёЄрты ыю ёюсющ
$m_{1}\times m_{2}$-ьрЄЁшЎє. ╧єёЄ№, ъръ ш т√°х, $\Phi_{}$ ---
фшёъЁхЄэюх яЁхюсЁрчютрэшх ╘єЁ№х, юяЁхфхыхээюх ЁртхэёЄтюь
(\ref{KorSlo.1.4}). ╨рёёьюЄЁшь {\it фшёъЁхЄэ√щ\/} ╧─╬
$f\Phi_{}\V\Phi_{}^{*}g$, фхщёЄтє■∙шщ шч
$L_{2}(\T^{d},\mathbb{C}^{m_{2}})$ т
$L_{2}(\T^{d},\mathbb{C}^{m_{1}})$. ═шцх фрэ√ эхюсїюфшь√х юЎхэъш ш
рёшьяЄюЄшъш ёшэуєы Ёэ√ї ўшёхы юяхЁрЄюЁр $f\Phi_{}\V\Phi_{}^{*}g$.
╥ръшх тюяЁюё√ шчєўрышё№ тю ьэюушї ЁрсюЄрї фы  юяхЁрЄюЁр
$\Tilde{f}\mathcal{F} \Tilde{\V}\mathcal{F}^{*}\Tilde{g}$, уфх
$\Tilde{f}(k)$, $\Tilde{\V}(\xi)$, $\Tilde{g}(k)$, $k\in\T^{d}$,
$\xi\in\rd$, --- ьрЄЁшЎ√-ЇєэъЎшш, ёюуырёютрээ√ї яюЁ фъют,
$\mathcal{F}$
--- яЁхюсЁрчютрэшх ╘єЁ№х т
$L_{2}(\rd)$. ╠√ ёыхфєхь шфх ь ЁрсюЄ 
╠.~╪.~┴шЁьрэр, ╠.~╟.~╤юыюь ър \cite{KorSlo.BirmanSolomyak90}, 
\cite{KorSlo.BirmanSolomyak77}, \cite{KorSlo.BirmanSolomyak79}, р
Єръ цх ЁрсюЄ√ \cite{KorSlo.Sloushch2014}. ╧Ёш юЎхэъх ёшэуєы Ёэ√ї
ўшёхы юяхЁрЄюЁр $f\Phi_{}\V\Phi_{}^{*}g$ ь√ яюы№чєхьё  яЁшхьюь шч
\cite{KorSlo.RozenblumSolomyak2009}.

\setcounter{nnn}{0} \pr {\bf ═хюсїюфшь√х ётхфхэш  ю ъюьяръЄэ√ї
юяхЁрЄюЁрї.}
─ы  яЁюшчтюы№эюую ъюьяръЄэюую юяхЁрЄюЁр $\A$, фхщёЄтє■∙хую шч
ушы№схЁЄютр яЁюёЄЁрэёЄтр $\H_{1}$ т ушы№схЁЄютю яЁюёЄЁрэёЄтю
$\H_{2}$, юсючэрўшь ўхЁхч $s_{m}(\A)$, $m\in\NN$, ёшэуєы Ёэ√х
ўшёыр юяхЁрЄюЁр $\A$ (Є.х. яюёыхфютрЄхы№э√х ёюсёЄтхээ√х чэрўхэш 
юяхЁрЄюЁр $|\A|:=(\A^{*}\A)^{1/2}$). ╬яЁхфхышь
$n(s,\A):=\#\{m\in\NN:s_{m}(\A)>s\}$ --- ЇєэъЎш■ ЁрёяЁхфхыхэш 
ёшэуєы Ёэ√ї ўшёхы юяхЁрЄюЁр $\A$. ─ы  ёрьюёюяЁ цхээюую ъюьяръЄэюую
юяхЁрЄюЁр $\A$ т ушы№схЁЄютюь яЁюёЄЁрэёЄтх $\H$ тхЁэю ЁртхэёЄтю
$n_{\pm}(s,\A)=n(s,\A_{\pm})$, уфх $\A_{\pm}:=(|\A|\pm\A)/2$.
╤яЁртхфыштю эхёъюы№ъю юўхтшфэ√ї ёююЄэю°хэшщ:
\begin{gather*}
n(s,\A^{*}\A)=n(\sqrt{s},\A)=n(\sqrt{s},\A^{*})=n(s,\A\A^{*}),\ \
s>0; \\
n(s,\A)=n_{+}(s,\A)+n_{-}(s,\A),\ \ s>0,\ \ \A=\A^{*};\\
n(s,\A)=n_{+}(s,\A),\ \ s>0,\ \ \A\ge 0.
\end{gather*}
─ы  ЇєэъЎшш ЁрёяЁхфхыхэш  $s$-ўшёхы ёєьь√ фтєї ъюьяръЄэ√ї
юяхЁрЄюЁют $\A$ ш $\mathbb{B}$ ёяЁртхфыштю эхЁртхэёЄтю ╘рэ№-╓■ 
(ёь., эряЁшьхЁ, \cite{KorSlo.BSU}, \S 11.1, я.3)
\begin{equation}\label{KorSlo.2.1.1}
n(s+t,\A+\mathbb{B})\le n(s,\A)+n(t,\mathbb{B}),\ \ s,t>0.
\end{equation}
╬ЄьхЄшь шёяюы№чєхьюх т фры№эхщ°хь <<трЁшрЎшюээюх>> ётющёЄтю
ЇєэъЎшш ЁрёяЁхфхыхэш  $n_{\pm}(s,\A)$ ёрьюёюяЁ цхээюую ъюьяръЄэюую
юяхЁрЄюЁр $\A$ т ушы№схЁЄютюь яЁюёЄЁрэёЄтх $\H$ (ёь., эряЁшьхЁ,
\cite{KorSlo.BSU}, \S 10.2, я.2):
\begin{equation}\label{KorSlo.2.1}
n_{\pm}(s,\A)=\sup\,\{\mathrm{dim}\,\mathscr{F},\
\mathscr{F}\subset\H : \pm(\A u,u)>s\|u\|^{2},\ \forall
u\in\mathscr{F}\setminus\{0\}\}.
\end{equation}

╩ырёё $\SS_{p,\infty}(\H_{1},\H_{2})$ (ёь., эряЁшьхЁ,
\cite{KorSlo.BirmanSolomyakSnumbers77}) т√фхы хЄё  єёыютшхь
ъюэхўэюёЄш ЇєэъЎшюэрыр
$$
\|\A\|_{\SS_{p,\infty}}:=\sup_{s>0}sn^{1/p}(s,\A).
$$
 ╧ЁюёЄЁрэёЄтю
$\SS_{p,\infty}$ яюыэю юЄэюёшЄхы№эю ътрчшэюЁь√
$\|\cdot\|_{\SS_{p,\infty}}$, тююс∙х уютюЁ , эхёхярЁрсхы№эю ш
ёюфхЁцшЄ ёхярЁрсхы№эюх яюфяЁюёЄЁрэёЄтю
$$\SS_{p,\infty}^{0}:=\{\A\in\SS_{p,\infty}:n(s,\A)=o(s^{-p}),s\to+0\},$$
т ъюЄюЁюь яыюЄэю ьэюцхёЄтю юяхЁрЄюЁют ъюэхўэюую Ёрэур. ═р
яЁюёЄЁрэёЄтх $\SS_{p,\infty}$ юяЁхфхыхэ√ эхяЁхЁ√тэ√х (ёь.,
эряЁшьхЁ, \cite{KorSlo.BS83}) ЇєэъЎшюэры√
\begin{equation*}
\mathfrak{D}_{p}(\A):=\limsup_{s\to 0}s^{p}n(s,\A);\ \ \ \ \ \
\mathfrak{d}_{p}(\A):=\liminf_{s\to 0}s^{p}n(s,\A).
\end{equation*}
╬ЄьхЄшь, ўЄю ЁртхэёЄтю $\mathfrak{D}_{p}(\A)=\mathfrak{d}_{p}(\A)$
ючэрўрхЄ ёяЁртхфыштюёЄ№ рёшьяЄюЄшъш $$n(s,\A)\sim
s^{-p}\mathfrak{D}_{p}(\A),\ \ s\to+0.$$

═р ьэюцхёЄтх ёрьюёюяЁ цхээ√ї юяхЁрЄюЁют шч $\SS_{p,\infty}$
эхяЁхЁ√тэ√ ЇєэъЎшюэры√
$\mathfrak{D}_{p}^{\pm}(\A):=\limsup_{s\to+0}s^{p}n_{\pm}(s,\A)$,
$\mathfrak{d}_{p}^{\pm}(\A):=\liminf_{s\to+0}s^{p}n_{\pm}(s,\A)$.
┬ эхъюЄюЁ√ї ЇюЁьєышЁютърї сєфхЄ шёяюы№чютрЄ№ё  юсючэрўхэшх
$D_{p}(\A)$ фы  ы■сюую шч Їєэъ\-Ўш\-ю\-эр\-ыют
$\mathfrak{D}_{p}(\A)$, $\mathfrak{D}_{p}^{\pm}(\A)$,
$\mathfrak{d}_{p}(\A)$, $\mathfrak{d}_{p}^{\pm}(\A)$. ╧Ёш ¤Єюь
ёыєўрщ
$D_{p}(\A)=\mathfrak{D}_{p}^{\pm}(\A),\mathfrak{d}_{p}^{\pm}(\A)$
ртЄюьрЄшўхёъш яюфЁрчєьхтрхЄ ЁртхэёЄтю $\A=\A^{*}$. ═хЄЁєфэю
тшфхЄ№, ўЄю фы  юяхЁрЄюЁр $\A\in\SS_{p,\infty}$ тъы■ўхэшх
$\A\in\SS_{p,\infty}^{0}$ ¤ътштрыхэЄэю ЁртхэёЄтє
$\mathfrak{D}_{p}(\A)=0$. ╬ЄьхЄшь шёяюы№чєхьюх т фры№эхщ°хь
ётющёЄтю (ёь., эряЁшьхЁ, \cite{KorSlo.BS83})
\begin{equation*}
D_{p}(\A+\mathbb{K})=D_{p}(\A),\ \
D=\mathfrak{D},\mathfrak{D}^{\pm},\mathfrak{d},\mathfrak{d}^{\pm},\
\ \A\in\SS_{p,\infty},\ \ \mathbb{K}\in\SS_{p,\infty}^{0}.
\end{equation*}
═рь яюэрфюсшЄё  ёыхфє■∙хх єЄтхЁцфхэшх (ёь., эряЁшьхЁ,
\cite{KorSlo.BSU}, \S 11.6, я.3).
\begin{PPP}\label{KorSloPropositionWOGelder}
┼ёыш т√яюыэхэ√ єёыютш  $\A\in\SS_{p,\infty}$,
$\mathbb{B}\in\SS_{q,\infty}$, Єю тхЁэю тъы■ўхэшх
$$\A\mathbb{B}\in\SS_{r,\infty},\ \ \text{уфх}\ \ r^{-1}=p^{-1}+q^{-1};$$ яЁш
¤Єюь ёяЁртхфыштю эхЁртхэёЄтю
\begin{equation*}
\|\A\mathbb{B}\|_{\SS_{r,\infty}}\le
C(p,q)\|\A\|_{\SS_{p,\infty}}\|\mathbb{B}\|_{\SS_{q,\infty}}.
\end{equation*}
┼ёыш фюяюыэшЄхы№эю $\A\in\SS_{p,\infty}^{0}$ шыш
$\mathbb{B}\in\SS_{q,\infty}^{0}$, Єю
$\A\mathbb{B}\in\SS_{r,\infty}^{0}$.
\end{PPP}
\pr{\bf ╬Ўхэъш Єшяр ╓тшъхы  фы  юърщьыхээюую фшёъЁхЄэюую
яЁхюсЁрчютрэш  ╘єЁ№х.}
╧єёЄ№ $(\X,d\rho)$ ш $(\Y,d\tau)$ --- фтр ёхярЁрсхы№э√ї шчьхЁшь√ї
яЁюёЄЁрэёЄтр ё $\sigma$-ъюэхўэ√ьш ьхЁрьш. ╧Ёхфяюыюцшь
$T:L_{2}(\Y,d\tau)\to L_{2}(\X,d\rho)$ --- ышэхщэ√щ юуЁрэшўхээ√щ
шэЄхуЁры№э√щ юяхЁрЄюЁ ё  фЁюь $t(\cdot,\cdot)\in
L_{\infty}(\X\times\Y,d\rho\times d\tau)$. ═рё шэЄхЁхёє■Є єёыютш 
юуЁрэшўхээюёЄш, ъюьяръЄэюёЄш, р Єръцх юЎхэъш ёшэуєы Ёэ√ї ўшёхы
юяхЁрЄюЁр $fTg$ яЁш яюфїюф ∙шї шчьхЁшь√ї тхёрї
$f:\X\to\mathbb{C}$, $g:\Y\to\mathbb{C}$. ╬ЄтхЄ эр ¤Єш тюяЁюё√
фрхЄ юсюс∙хээр  {\it юЎхэър ╓тшъхы \/} (ёь.
\cite{KorSlo.BirmanSolomyak90}):
\begin{Th}\label{KorSloTheoremCwikelBirmanSolomyak}
╧єёЄ№  $(f, g)\in L_{p,\infty}(\X,d\rho)\times L_{p}(\Y,d\tau)$ яЁш
 $p>2$. ╥юуфр ёяЁртхфышт√ тъы■ўхэшх $fTg\in\SS_{p,\infty}$ ш юЎхэър
\begin{equation*}
\|fTg\|_{\SS_{p,\infty}}\le
C(p)\|T\|^{1-\frac{2}{p}}\|t\|_{L_{\infty}}^{\frac{2}{p}}\|f\|_{L_{p,\infty}}\|g\|_{L_{p}}.
\end{equation*}
┼ёыш фюяюыэшЄхы№эю $\rho_{f}(s)=o(s^{-p})$, $s\to+0$, Єю
$fTg\in\SS_{p,\infty}^{0}$. ╨рчєьххЄё , єёыютш  эр тхёр $f$ ш $g$
ьюцэю яюьхэ Є№ ьхёЄрьш.
\end{Th}
─рыхх юЄьхЄшь, ўЄю фы  $m_{1}\times m_{2}$-ьрЄЁшўэюую яюЄхэЎшрыр
$\W\in \ell_{\infty}(\zd)$, єс√тр■∙хую эр схёъюэхўэюёЄш, юяхЁрЄюЁ
$[\W(n)]$ ъюьяръЄхэ. ┴юыхх Єюую, ёяЁртхфыштю ёыхфє■∙хх
єЄтхЁцфхэшх.

\begin{Nz}\label{KorSlo.Note2.3}
┼ёыш $\W\in \ell_{p,\infty}(\zd)$, Єю $[\W(n)]\in\SS_{p,\infty}$,
$\|[\W(n)]\|_{\SS_{p,\infty}}\le
C(m_{1},m_{2})\|\W\|_{\ell_{p,\infty}}$.

 ┼ёыш $\W\in
\ell_{p,\infty}^{\,0}(\zd)$, Єю $[\W(n)]\in\SS_{p,\infty}^{0}$.
\end{Nz}

╨рёёьюЄЁшь яЁ ьюєуюы№э√х ьрЄЁшЎ√-ЇєэъЎшш ёюуырёютрээ√ї яюЁ фъют
$f(k)$, $k\in\T^{d}$, $\W(n)$, $n\in\zd$. ╤ыхфє■∙хх єЄтхЁцфхэшх
фрхЄ єёыютш  ъюьяръЄэюёЄш, р Єръцх юЎхэъш ёшэуєы Ёэ√ї ўшёхы
юяхЁрЄюЁр $f\Phi_{}\W$ (юЎхэъє Єшяр ╓тшъхы ).
\begin{PPP}\label{KorSloPropositionCwikelTypeEstimate}
╧єёЄ№ тхЁэ√ тъы■ўхэш  $f\in L_{q}(\T^{d})$, $\W\in
\ell_{p,\infty}(\zd)$, уфх $p\in(0,+\infty)$ ш
\begin{enumerate}
\item[р)] $q=p$ яЁш $p\in(2,+\infty)$; \item[с)] $q=2$ яЁш
$p\in(0,2)$; \item[т)] $q>2$ яЁш $p=2$.
\end{enumerate}
╥юуфр ёяЁртхфышт√ тъы■ўхэшх $f\Phi_{}\W\in\SS_{p,\infty}$ ш юЎхэър
\begin{equation}\label{KorSlo.2.2}
\|f\Phi_{}\W\|_{\SS_{p,\infty}}\le
C\|f\|_{L_{q}}\|\W\|_{\ell_{p,\infty}}.
\end{equation}
╩юэёЄрэЄр $C=C(p,q,d)$ т юЎхэъх {\rm(\ref{KorSlo.2.2})} чртшёшЄ
Єръцх юЄ яюЁ фъют ьрЄЁшЎ $f$ ш $\W$. ┼ёыш фюяюыэшЄхы№эю $\W\in
\ell_{p,\infty}^{\,0}(\zd)$, Єю $f\Phi_{}\W\in\SS_{p,\infty}^{0}$.
\end{PPP}
\begin{proof}[─юърчрЄхы№ёЄтю]
─юёЄрЄюўэю ЁрёёьюЄЁхЄ№ ёыєўрщ ёъры Ёэ√ї ЇєэъЎшщ $f(k)$, $\W(n)$.
\par\noindent{\rm 1)} ╧Ёш $p\in(2,+\infty)$ ЄЁхсєхь√х єЄтхЁцфхэш 
т√Єхър■Є шч ЄхюЁхь√ \ref{KorSloTheoremCwikelBirmanSolomyak}.
\par\noindent{\rm 2)}
┬ ёыєўрх $p\in(0,2)$ фы  яЁютхЁъш тъы■ўхэш 
$f\Phi\W\in\SS_{p,\infty}$ ш юЎхэъш {\rm(\ref{KorSlo.2.2})}
шёяюы№чє■Єё  рЁуєьхэЄ√ ЁрсюЄ√ \cite{KorSlo.RozenblumSolomyak2009}.
╚ьхээю, т ёшыє тъы■ўхэш  $\ell_{p,\infty}(\zd)\subset
\ell_{2}(\zd)$, $p\in(0,2)$, юяхЁрЄюЁ $f\Phi\W$
--- ъырёёр ├шы№схЁЄр-╪ьшфЄр. ─юёЄрЄюўэю яЁютхЁ Є№ тъы■ўхэшх $f\Phi\W\in\SS_{p,\infty}$ ш юЎхэъє (\ref{KorSlo.2.2}) т яЁхфяюыюцхэшш $\|f\|_{L_{2}}=1$.
╨рчюс№хь ЇєэъшЎ■ $\W$ эр фтх ўрёЄш $\W=\W^{s}+\W_{s}$, $s>0$;
чфхё№ $\W_{s}(n)=\W(n)$, хёыш $|\W(n)|\le s$, ш $\W_{s}(n)=0$,
хёыш $|\W(n)|>s$. ╬сючэрўшт ўхЁхч $d\mu$ ёўшЄр■∙є■ ьхЁє эр $\zd$,
юЄьхЄшь ¤ыхьхэЄрЁэ√х ётющёЄтр
\begin{equation}\label{KorSlo.2.3}
\left\{
\begin{array}{l}
\mu_{\,\W_{s}}(\sigma)=0,\ \ \sigma\ge s,\\
\mu_{\,\W_{s}}(\sigma)\le\mu_{\,\W}(\sigma),\ \ \sigma<s;
\end{array}
\right.\ \ \mathrm{rank}\,\W^{s}=\mu_{\,\W}(s),\ \ s>0.
\end{equation}
┬ ёшыє (\ref{KorSlo.2.1.1}) ёяЁртхфыштю эхЁртхэёЄтю
\begin{equation}\label{KorSlo.2.4}
\textstyle
n(s,f\Phi\W)\le
n(\frac{s}{2},f\Phi\W^{s})+n(\frac{s}{2},f\Phi\W_{s}),\ \ s>0.
\end{equation}
╧юёъюы№ъє $\W^{s}$ --- юяхЁрЄюЁ ъюэхўэюую Ёрэур (ёь.
(\ref{KorSlo.2.3})), яхЁтюх ёырурхьюх т яЁртющ ўрёЄш
(\ref{KorSlo.2.4}) фюяєёърхЄ юЎхэъє
\begin{equation}\label{KorSlo.2.5}
\textstyle
n(\frac{s}{2},f\Phi\W^{s})\le\mu_{\,\W}(s)\le\|\W\|_{\ell_{p,\infty}}^{p}s^{-p},\
\ s>0.
\end{equation}
╧юёъюы№ъє $\W_{s}\in\ell_{2}$, тЄюЁюх ёырурхьюх т яЁртющ ўрёЄш
(\ref{KorSlo.2.4}) єфютыхЄтюЁ хЄ эхЁртхэёЄтє
\begin{equation}\label{KorSlo.2.6}
\textstyle
n(\frac{s}{2},f\Phi\W_{s})\le
4s^{-2}\|f\Phi\W_{s}\|_{\SS_{2}}^{2}=4s^{-2}(2\pi)^{-d}\|\W_{s}\|_{\ell_{2}}^{2}
\end{equation}
╬ёЄрхЄё  чрьхЄшЄ№, ўЄю шч (\ref{KorSlo.2.3}) т√Єхър■Є ёююЄэю°хэш 
\begin{equation}\label{KorSlo.2.7}
\|\W_{s}\|_{\ell_{2}}^{2}=2\int_{0}^{+\infty}\sigma\mu_{\,\W_{s}}(\sigma)d\sigma\le
2\int_{0}^{s}\sigma\mu_{\,\W}(\sigma)d\sigma\le\frac{2}{2-p}\|\W\|_{\ell_{p,\infty}}^{p}s^{2-p}.
\end{equation}
╥хяхЁ№ тъы■ўхэшх $f\Phi\W\in\SS_{p,\infty}$ ш юЎхэър
(\ref{KorSlo.2.2}) т√Єхър■Є шч (\ref{KorSlo.2.4}) --
(\ref{KorSlo.2.7}).
\par\noindent{\rm 3)}
═ръюэхЎ, т ёыєўрх $p=2$ тюёяюы№чєхьё  ЁртхэёЄтюь
\begin{equation*}
\|f\Phi_{}\W\|_{2,\infty}^{2}=\|f\Phi_{}|\W|^{2\theta_{1}}|\W|^{2\theta_{2}}\Phi_{}^{*}\bar{f}\|_{1,\infty}^{},\
\ \theta_{1}+\theta_{2}=1, \theta_{2}=1/q.
\end{equation*}
╬ёЄрхЄё  яЁшьхэшЄ№ ъ юяхЁрЄюЁрь $f\Phi_{}|\W|^{2\theta_{1}}$,
$|\W|^{2\theta_{2}}\Phi_{}^{*}\bar{f}$ юЎхэъє (\ref{KorSlo.2.2})
фы  ёыєўрхт $p=1/\theta_{1}\in(0,2)$ ш $p=1/\theta_{2}=q>2$,
ёююЄтхЄёЄтхээю, ш тюёяюы№чютрЄ№ё  яЁхфыюцхэшхь
\ref{KorSloPropositionWOGelder}.
\end{proof}
\begin{PPP}\label{KorSloPropositionDiscreteCwikelEstimate}
╧єёЄ№ яЁш эхъюЄюЁюь $p\in(2,+\infty)$ т√яюыэхэ√ єёыютш  $f\in
L_{p,\infty}(\T^{d})$, $\W\in\ell_{p}(\zd)$. ╥юуфр ёяЁртхфышт√
тъы■ўхэшх $f\Phi_{}\W\in\SS_{p,\infty}^{0}$ ш юЎхэър
\begin{equation}\label{KorSlo.2.8}
\|f\Phi_{}\W\|_{\SS_{p,\infty}}\le
C(p,d)\|f\|_{L_{p,\infty}}\|\W\|_{\ell_{p}}.
\end{equation}
╩юэёЄрэЄр $C(p,d)$ т юЎхэъх {\rm(\ref{KorSlo.2.8})} чртшёшЄ юЄ
яюЁ фъют ьрЄЁшЎ $f$ ш $\W$.
\end{PPP}
\begin{proof}
╧Ёш яЁютхЁъх яЁхфыюцхэш 
\ref{KorSloPropositionDiscreteCwikelEstimate} фюёЄрЄюўэю, ъръ ш
т√°х, ЁрёёьюЄЁхЄ№ ёыєўрщ ёъры Ёэ√ї ЇєэъЎшщ $f(k)$, $k\in\T^{d}$,
$\W(n)$, $n\in\zd$. ЄЁхсєхь√х єЄтхЁцфхэш  т√Єхър■Є шч ЄхюЁхь√
\ref{KorSloTheoremCwikelBirmanSolomyak} ш юЎхэъш
$\mathrm{mes}\{k\in\T^{d}:|f(k)|>s\}=o(s^{-p}),\ \ s\to+0$.
\end{proof}
╚ч яЁхфыюцхэшщ \ref{KorSloPropositionCwikelTypeEstimate} ш
\ref{KorSlo.Note1.1} т√ЄхърхЄ єЄтхЁцфхэшх.
\begin{Con}\label{KorSloCorollaryCwikelTypeEstimate}
╧єёЄ№ т√яюыэхэ√ єёыютш  $\W(n)=O(|n|^{-d/p})$, $p>0$, $f\in
L_{q}(\T^{d})$, уфх $q$ Єръюх цх, ъръ т яЁхфыюцхэшш
{\rm\ref{KorSloPropositionCwikelTypeEstimate}}.
╥юуфр ёяЁртхфыштю тъы■ўхэшх
$f\Phi_{}\W\Phi_{}^{*}\in\SS_{p,\infty}$.
\end{Con}
╧Ёштхфхь ётющёЄтю яЁшсышцхээюую ъюььєЄшЁютрэш  юяхЁрЄюЁют $f$ ш
$\Phi_{}\W\Phi_{}^{*}$ яЁш яюфїюф ∙шї $f$ ш $\W$ (сышчъшщ, эю
Єхїэшўхёъш сюыхх ёыюцэ√щ Ёхчєы№ЄрЄ фы  эхяЁхЁ√тэюую ёыєўр  с√ы
фюърчрэ т \cite{KorSlo.Sloushch2014}).
\begin{PPP}\label{KorSloPropositionApproxComm}
╧єёЄ№ $\W(n)$, $n\in\zd$, ш $f(k)$, $k\in\T^{d}$, ътрфЁрЄэ√х
ьрЄЁшЎ√ юфшэръют√ї яюЁ фъют; яєёЄ№ яЁш эхъюЄюЁюь $p\in(0,+\infty)$
ёяЁртхфыштр рёшьяЄюЄшър
\begin{equation*}
\textstyle
\W(n)=|n|^{-d/p}\left(\omega(\frac{n}{|n|})\mathbbm{1}+o(1)\right),\
\ |n|\to+\infty.
\end{equation*}
╟фхё№ $\mathbbm{1}$ --- хфшэшўэр  ьрЄЁшЎр, $\omega\in
C(\mathbb{S}^{d-1})$ --- ёъры Ёэр  ЇєэъЎш . ╧єёЄ№ $f\in
L_{q}(\T^{d})$, уфх $q$ Єръюх цх, ъръ т яЁхфыюцхэшш
{\rm\ref{KorSloPropositionCwikelTypeEstimate}}.
╥юуфр тхЁэю тъы■ўхэшх
\begin{equation}\label{KorSlo.2.9}
f\Phi_{}\W\Phi_{}^{*}-\Phi_{}\W\Phi_{}^{*}f\in\SS_{p,\infty}^{0}.
\end{equation}
\end{PPP}
\begin{proof}
╤юуырёэю яЁхфыюцхэш■ \ref{KorSlo.Note1.1} ш яЁхфыюцхэш■
\ref{KorSloPropositionCwikelTypeEstimate} тъы■ўхэшх
(\ref{KorSlo.2.9}) фюёЄрЄюўэю яЁютхЁ Є№ т яЁхфяюыюцхэшш $\W(0)=0$,
$\W(n)=\omega(\frac{n}{|n|})\mathbbm{1}|n|^{-d/p}$, $n\not=0$. ╧Ёш
¤Єюь ьрЄЁшЎр $\W(n)$, $n\in\zd$, яЁюяюЁЎшюэры№эр хфшэшўэющ, р
яюЄюьє фюёЄрЄюўэю ЁрёёьюЄЁхЄ№ ёыєўрщ ёъры Ёэ√ї ЇєэъЎшщ $\W(n)$ ш
$f(k)$. ╥хяхЁ№, т ёшыє яЁхфыюцхэш  \ref{KorSlo.Note1.1} ш
яЁхфыюцхэш  \ref{KorSloPropositionCwikelTypeEstimate} тъы■ўхэшх
(\ref{KorSlo.2.9}) ьюцэю яЁютхЁ Є№, яЁхфяюырур  $\omega\in
C^{\infty}(\mathbb{S}^{d-1})$,
$f(k)=(2\pi)^{-d/2}\sum\limits_{|n|\le N}f_{n}e^{-ink}$,
$k\in\T^{d}$, $N\in\mathbb{N}$.

╧юёъюы№ъє юяхЁрЄюЁ $\Phi_{}$ єэшЄрЁхэ, тъы■ўхэшх
(\ref{KorSlo.2.9}) ¤ътштрыхэЄэю ёююЄэю°хэш■
\begin{equation}\label{KorSlo.2.10}
\Phi_{}^{*}f\Phi_{}\W-\W\Phi_{}^{*}f\Phi_{}\in\SS_{p,\infty}^{0}.
\end{equation}
─ююяЁхфхышь $f_{n}=0$ яЁш $|n|>N$. ╬яхЁрЄюЁ $\Phi_{}^{*}f\Phi_{}$
--- шэЄхуЁры№э√щ (ёєььрЄюЁэ√щ) юяхЁрЄюЁ ё  фЁюь
$(2\pi)^{-d/2}f_{n-m}$, Є.х.
\begin{equation*}
\Phi_{}^{*}f\Phi_{}u(n)=(2\pi)^{-d/2}\sum_{m\in\zd:|n-m|\le
N}f_{n-m}u(m),\ \ u\in \ell_{2}(\zd).
\end{equation*}
╤ыхфютрЄхы№эю, юяхЁрЄюЁ
$\Phi_{}^{*}f\Phi_{}\W-\W\Phi_{}^{*}f\Phi_{}$ --- шэЄхуЁры№э√щ
юяхЁрЄюЁ ё  фЁюь $$(2\pi)^{-d/2}f_{n-m}(\W(m)-\W(n)),\ \
m,n\in\zd;$$ юЄё■фр эхЄЁєфэю т√тхёЄш ЁртхэёЄтю
\begin{equation}\label{KorSlo.2.11}
\Phi_{}^{*}f\Phi_{}\W-\W\Phi_{}^{*}f\Phi_{}=(2\pi)^{-d/2}\sum_{|t|\le
N}f_{t}S_{-t}[\W(n)-\W(n+t)].
\end{equation}
╟фхё№ $S_{a}u(x)=u(x+a)$, $a\in\zd$, --- єэшЄрЁэ√щ юяхЁрЄюЁ ёфтшур
т $\ell_{2}(\zd)$. ╬ёЄрхЄё  чрьхЄшЄ№, ўЄю ёяЁртхфыштр юЎхэър
$\W(n)-\W(n+t)=o(|n|^{-d/p})$, $|n|\to+\infty$, $t\in\zd$.
╧юёыхфэхх ёююЄэю°хэшх тьхёЄх ё (\ref{KorSlo.2.11}) яЁхфыюцхэшхь
\ref{KorSlo.Note1.1} ш чрьхўрэшхь \ref{KorSlo.Note2.3} яЁштюфшЄ ъ
(\ref{KorSlo.2.10}).
\end{proof}
\pr {\bf └ёшьяЄюЄшър ёшэуєы Ёэ√ї ўшёхы фшёъЁхЄэюую ╧─╬
юЄЁшЎрЄхы№эюую яюЁ фър.} ═шцх $f(k)$, $\V(n)$, $g(k)$,
$k\in\T^{d}$, $n\in\zd$,  --- яЁ ьюєуюы№э√х ьрЄЁшЎ√-ЇєэъЎшш
ёюуырёютрээ√ї яюЁ фъют.
\begin{Th}\label{KorSloTheoremPDOEstimateBirmanSolomyak}
╧єёЄ№ т√яюыэхэ√ єёыютш  $\V\in \ell_{p,\infty}(\zd)$, $f\in
L_{q_{1}}(\T^{d})$, $g\in L_{q_{2}}(\T^{d})$; чфхё№
$p\in(0,+\infty)$ ш
\begin{enumerate}
\item[р)] хёыш $p\in(1,+\infty)$, Єю $q_{1},q_{2}\in(2,+\infty]$,
$1/q_{1}+1/q_{2}=1/p$;
\item[с)] хёыш $p\in(0,1)$, Єю $q_{1}=q_{2}=2$,
\item[т)] хёыш $p=1$, Єю ышсю $q_{1}=2,\ \ q_{2}>2$, ышсю
$q_{2}=2,\ \ q_{1}>2$.
\end{enumerate}
╥юуфр ёяЁртхфышт√ тъы■ўхэшх
$f\Phi_{}\V\Phi_{}^{*}g\in\SS_{p,\infty}$ ш юЎхэър
\begin{equation}\label{KorSlo.2.12}
\|f\Phi_{}\V\Phi_{}^{*}g\|_{\SS_{p,\infty}}\le
C\|f\|_{L_{q_{1}}}\|g\|_{L_{q_{2}}}\|\V\|_{\ell_{p,\infty}}.
\end{equation}
╩юэёЄрэЄр $C$ т {\rm(\ref{KorSlo.2.12})} чртшёшЄ юЄ
$p,q_{1},q_{2},d$ ш юЄ яюЁ фър ьрЄЁшЎ $f$, $\V$ ш $g$. ┼ёыш
фюяюыэшЄхы№эю $\V\in \ell_{p,\infty}^{\,0}(\zd)$, Єю
$f\Phi_{}\V\Phi_{}^{*}g\in\SS_{p,\infty}^{0}$.
\end{Th}
\begin{proof}
─юёЄрЄюўэю ЁрёёьюЄЁхЄ№ ёъры Ёэ√щ ёыєўрщ. ┬юёяюы№чєхьё  Ёрчыюцхэшхь
\begin{equation*}
f\Phi_{}\V\Phi_{}^{*}g=f\Phi_{}|\V|^{\theta_{1}}\frac{\V}{|\V|}|\V|^{\theta_{2}}\Phi_{}^{*}g,\
\ \theta_{1}+\theta_{2}=1.
\end{equation*}
╧Ёш $p\in(1,+\infty)$ яюыюцшь $\theta_{i}=p/q_{i}$, $i=1,2$; яЁш
$p\in(0,1)$ т√схЁхь $\theta_{1}=\theta_{2}=1/2$; т ёыєўрх $p=1$
(ш, эряЁшьхЁ, $q_{1}=2$, $q_{2}>2$) т√схЁхь $\theta_{2}=p/q_{2}$,
$\theta_{1}=1-\theta_{2}$. ─рыхх, яЁшьхэшь ъ юяхЁрЄюЁрь
$f\Phi_{}|\V|^{\theta_{1}}$ ш $|\V|^{\theta_{2}}\Phi_{}^{*}g$
яЁхфыюцхэшх \ref{KorSloPropositionCwikelTypeEstimate} ш, єўшЄ√тр 
яЁхфыюцхэшх \ref{KorSloPropositionWOGelder}, яюыєўшь ЄЁхсєхь√х
єЄтхЁцфхэш .
\end{proof}
╚ч яЁхфыюцхэш  \ref{KorSloPropositionApproxComm} ш ЄхюЁхь√
\ref{KorSloTheoremPDOEstimateBirmanSolomyak} т√ЄхърхЄ ёыхфє■∙хх
єЄтхЁцфхэшх.
\begin{Con}\label{KorSlo.Corollary2.8} ╧єёЄ№ т√яюыэхэ√ єёыютш  ЄхюЁхь√
{\rm\ref{KorSloTheoremPDOEstimateBirmanSolomyak}} ш фюяюыэшЄхы№эю
\begin{equation}\label{KorSlo.2.13}
\textstyle
\V(n)=|n|^{-d/p}\left(v(\frac{n}{|n|})+o(1)\right),\ \ v\in
C(\mathbb{S}^{d-1}),
\end{equation}
ш эюёшЄхыш ЇєэъЎшщ $f$ ш $g$ эх яхЁхёхър■Єё . ╥юуфр ёяЁртхфыштю
тъы■ўхэшх
\begin{equation}\label{KorSlo.2.14}
f\Phi_{}\V\Phi_{}^{*}g\in\SS_{p,\infty}^{0}.
\end{equation}
\end{Con}
\begin{proof}
─юёЄрЄюўэю ЁрёёьюЄЁхЄ№ ёыєўрщ ёъры Ёэ√ї ЇєэъЎшщ $f(k)$, $\V(n)$,
$g(k)$, $k\in\T^{d}$, $n\in\zd$. ─рыхх, яЁхфяюыюцшь фы 
юяЁхфхыхээюёЄш $q_{1}<\infty$. ╤юуырёэю ЄхюЁхьх
\ref{KorSloTheoremPDOEstimateBirmanSolomyak}, фюёЄрЄюўэю яЁютхЁ Є№
(\ref{KorSlo.2.14}) яЁш $f\in C_{0}^{\infty}(\T^{d})$, $g\in
L_{\infty}(\T^{d})$; т ¤Єшї яЁхфяюыюцхэш ї (\ref{KorSlo.2.14})
т√ЄхърхЄ шч яЁхфыюцхэш 
\ref{KorSloPropositionApproxComm}. 
\end{proof}
═шцх фы  яЁюшчтюы№эющ ьрЄЁшЎ√ $B$ яюыюцшь
$\Lambda_{p}(B):=\sum\limits_{k}s_{k}^{p}(B)$, $p\in(0,+\infty)$.
\begin{Th}\label{KorSloTheoremPDOAsymptoticBirmanSolomyak} ╧єёЄ№ т√яюыэхэ√ єёыютш  ЄхюЁхь√
{\rm\ref{KorSloTheoremPDOEstimateBirmanSolomyak}},
$q_{1},q_{2}<\infty$, ш ёяЁртхфыштю ёююЄэю°хэшх
{\rm(\ref{KorSlo.2.13})}. ╥юуфр шьххЄ ьхёЄю рёшьяЄюЄшър $s$-ўшёхы
\begin{equation}\label{KorSlo.2.15}
\mathfrak{D}_{p}(f\Phi_{}\V\Phi_{}^{*}g)=\mathfrak{d}_{p}(f\Phi_{}\V\Phi_{}^{*}g)
=\frac{1}{d(2\pi)^{d}}\int_{\T^{d}}dk\int_{\mathbb{S}^{d-1}}\Lambda_{p}(f(k)v(\theta)g(k))\,dS(\theta).
\end{equation}
\end{Th}
\begin{proof} ╧Ёш $d=1$ ЄхюЁхьр \ref{KorSloTheoremPDOAsymptoticBirmanSolomyak} с√ыр фюърчрэр т ЁрсюЄрї \cite{KorSlo.HLN2016},
\cite{KorSlo.LN2018}. ╟фхё№ ь√ яЁштхфхь фюърчрЄхы№ёЄтю фы 
яЁюшчтюы№эющ ЁрчьхЁэюёЄш. ╩ръ ш т√°х, ьюцэю юуЁрэшўшЄ№ё  ёыєўрхь
$\V(n)=v(\frac{n}{|n|})|n|^{-d/p}$, $n\not=0$, $\V(0)=0$. ┬ ёшыє
ЄхюЁхь√ \ref{KorSloTheoremPDOEstimateBirmanSolomyak}
рёшьяЄюЄшўхёъшх ъю¤ЇЇшЎшхэЄ√
$\mathfrak{D}_{p}(f\Phi_{}\V\Phi_{}^{*}g)$,
$\mathfrak{d}_{p}(f\Phi_{}\V\Phi_{}^{*}g)$ эхяЁхЁ√тэ√ юЄэюёшЄхы№эю
ъю¤ЇЇшЎшхэЄют $f$, $v$, $g$ т яЁюёЄЁрэёЄтх
$L_{q_{1}}(\T^{d})\times C(\mathbb{S}^{d-1})\times
L_{q_{2}}(\T^{d})$. ╤ фЁєующ ёЄюЁюэ√ яЁртр  ўрёЄ№ т
(\ref{KorSlo.2.15}) Єръцх эхяЁхЁ√тэр юЄэюёшЄхы№эю ъю¤ЇЇшЎшхэЄют
$f$, $v$, $g$ т яЁюёЄЁрэёЄтх $L_{q_{1}}(\T^{d})\times
C(\mathbb{S}^{d-1})\times L_{q_{2}}(\T^{d})$. ╤ыхфютрЄхы№эю,
фюёЄрЄюўэю яЁютхЁшЄ№ (\ref{KorSlo.2.15}), фы  $v\in
C^{\infty}(\mathbb{S}^{d-1})$ ш яЁюёЄ√ї ЇєэъЎшщ $f=\sum
f_{k}\mathbf{1}_{\mathscr{D}_{k}}$, $g=\sum
g_{k}\mathbf{1}_{\mathscr{D}_{k}}$, юЄтхўр■∙шї ъюэхўэюьє Ёрчсшхэш■
ъєср $\T^{d}$ эр яхЁхёхър■∙шхё  Єюы№ъю яю уЁрэшЎх ъєс√
$\{\mathscr{D}_{k}\}$. ┬ ¤Єюь ёыєўрх ёяЁртхфыштю ЁртхэёЄтю
\begin{equation}\label{KorSlo.2.16}
f\Phi_{}\V\Phi_{}^{*}g=\sum_{k,l}f_{k}\mathbf{1}_{\mathscr{D}_{k}}\Phi_{}\V\Phi_{}^{*}g_{l}\mathbf{1}_{\mathscr{D}_{l}}.
\end{equation}
╧Ёш $k\not=l$ ёырурхь√х т (\ref{KorSlo.2.16}) яЁшэрфыхцрЄ ъырёёє
$\SS_{p,\infty}^{0}$ (ёь. (\ref{KorSlo.2.14})). ╥ръшь юсЁрчюь,
т√яюыэхэю ЁртхэёЄтю
\begin{equation}\label{KorSlo.2.17}
\begin{aligned}
&
\mathcal{D}_{p}(f\Phi_{}\V\Phi_{}^{*}g)=\mathcal{D}_{p}(\sum_{k}f_{k}\mathbf{1}_{\mathscr{D}_{k}}
\Phi_{}\V\Phi_{}^{*}g_{k}\mathbf{1}_{\mathscr{D}_{k}})=
\sum_{k}\mathcal{D}_{p}(f_{k}\mathbf{1}_{\mathscr{D}_{k}}
\Phi_{}\V\Phi_{}^{*}g_{k}\mathbf{1}_{\mathscr{D}_{k}}),\\
& \mathcal{D}_{p}(T)=\lim_{s\to 0}s^{p}n(s,T).
\end{aligned}
\end{equation}
─ююяЁхфхышь ьрЄЁшЎє-ЇєэъЎш■ $\V\in \ell_{p,\infty}(\zd)$ тю тёхь
яЁюёЄЁрэёЄтх $\rd$ ЁртхэёЄтюь
$$\V(\xi)=\zeta(|\xi|)v(\xi/|\xi|)|\xi|^{-d/p},\ \ \xi\in\rd;$$
чфхё№ $\zeta(|\xi|)$ --- яюфїюф ∙р  уырфър  ёЁхчър т эєых.
╬сючэрўшь (ъръ ш т√°х) ўхЁхч ${\mathcal{F}}$ єэшЄрЁэюх
яЁхюсЁрчютрэшх ╘єЁ№х т $L_{2}(\rd)$:
\begin{equation*}
{\mathcal{F}}u(x)=(2\pi)^{-d/2}\int_{\rd}e^{-ixy}u(y)dy, \ \ u\in
L_{2}(\rd).
\end{equation*}
┬ ЁрсюЄх \cite{KorSlo.BirmanSolomyak79} (яя. 14--18) с√ыю фюърчрэю
ёююЄэю°хэшх
\begin{equation*}
\mathbf{1}_{\T^{d}}\Phi_{}\V\Phi_{}^{*}\mathbf{1}_{\T^{d}}-\mathbf{1}_{\T^{d}}{\mathcal{F}}{\V}{\mathcal{F}}^{*}\mathbf{1}_{\T^{d}}\in\SS_{p,\infty}^{0}.
\end{equation*}
╤ыхфютрЄхы№эю, ёяЁртхфыштю тъы■ўхэшх
\begin{equation}\label{KorSlo.2.18}
f_{k}\mathbf{1}_{\mathscr{D}_{k}}\Phi_{}\V\Phi_{}^{*}\mathbf{1}_{\mathscr{D}_{k}}g_{k}-
\mathbf{1}_{\mathscr{D}_{k}}{\mathcal{F}}f_{k}{\V}g_{k}{\mathcal{F}}^{*}\mathbf{1}_{\mathscr{D}_{k}}\in\SS_{p,\infty}^{0}.
\end{equation}
┬ ЁрсюЄрї \cite{KorSlo.BirmanSolomyak77},
\cite{KorSlo.BirmanSolomyak79} с√ыр эрщфхэр рёшьяЄюЄшър
\begin{equation}\label{KorSlo.2.19}
\mathcal{D}_{p}(\mathbf{1}_{\mathscr{D}_{k}}{\mathcal{F}}f_{k}{\V}g_{k}{\mathcal{F}}^{*}\mathbf{1}_{\mathscr{D}_{k}})=(2\pi)^{-d}d^{-1}\int_{\mathbb{S}^{d-1}}dS(\theta)
\int_{\mathscr{D}_{k}}\,{}\Lambda_{p}(f_{k}v(\theta)g_{k})\,dx.
\end{equation}
╤Ёртэштр  (\ref{KorSlo.2.19}), (\ref{KorSlo.2.18}) ш
(\ref{KorSlo.2.17}) яЁшїюфшь ъ (\ref{KorSlo.2.15}).
\end{proof}

\section{─юърчрЄхы№ёЄтю ЄхюЁхь \ref{KorSloTheoremRLC},
\ref{KorSloTheoremWRLC} ш \ref{KorSloTheoremAsymptotic}}

\setcounter{nnn}{0}  ─ы  юяЁхфхыхээюёЄш ь√ сєфхь яЁютхЁ Є№ юЎхэъш ш
рёшьяЄюЄшъш фы  тхышўшэ√ $N_{+}(\lambda,\tau)$,
$\lambda\in[\Lambda_{+},\Lambda_{-})$, $\tau>0$, т ёыєўрх {\it
тэєЄЁхээхщ\/} ыръєэ√ $(\Lambda_{+},\Lambda_{-})$ т ёяхъЄЁх
$\sigma(H)$, єфютыхЄтюЁ ■∙хщ яЁш эхъюЄюЁюь $N\in\{2,\dots,\nu\}$
єёыютш■ (\ref{KorSlo.1.7}).

\pr{\bf ╧ЁхфтрЁшЄхы№э√х чрьхўрэш .} 1) ╧Ёш єёыютш ї
(\ref{KorSlo.1.1}), (\ref{KorSlo.1.2}) ш (\ref{KorSlo.1.9}) (юэш
ёяЁртхфышт√ фы  ы■сющ шч ЄхюЁхь \ref{KorSloTheoremRLC},
\ref{KorSloTheoremWRLC}, \ref{KorSloTheoremAsymptotic}) ёўшЄр■∙р 
ЇєэъЎш  $N_{+}(\lambda,\tau)$ єфютыхЄтюЁ хЄ ЁртхэёЄтє
(\ref{KorSlo.1.14}$+$). ╬сючэрўшт фы  ъЁрЄъюёЄш
$E_{+}:=E_{H}(-\infty,\Lambda_{+}]$, шч (\ref{KorSlo.1.14}$+$) ш
трЁшрЎшюээюую ётющёЄтр (\ref{KorSlo.2.1}) т√тхфхь юЎхэъє
\begin{multline}\label{KorSlo.3.1}
N_{+}(\lambda,\tau)\le n(\tau^{-1},V^{1/2}E_{+}(\lambda
I-H)^{-1}E_{+}V^{1/2})=\\=n(\tau^{-1},\V^{1/2}\mathbb{E}_{+}(\lambda
I-\mathbb{H})^{-1}\mathbb{E}_{+}\V^{1/2}),\ \
\lambda\in(\Lambda_{+},\Lambda_{-}).
\end{multline}
╟фхё№ $\mathbb{H}=\mathbb{U}H\mathbb{U}^{*}$,
$\mathbb{V}=\mathbb{U}V\mathbb{U}^{*}$,
$\mathbb{E}_{+}=\mathbb{U}E_{+}\mathbb{U}^{*}=E_{\mathbb{H}}(-\infty,\Lambda_{+}]$,
єэшЄрЁэ√щ юяхЁрЄюЁ
$\mathbb{U}:\ell_{2}(X)\to\ell_{2}(\zd;\mathbb{C}^{\nu})$
юяЁхфхыхэ ЁртхэёЄтюь (\ref{KorSlo.1.3}).

╤юуырёэю яЁхфёЄртыхэш■ (\ref{KorSlo.1.5}) юяхЁрЄюЁ $\mathbb{H}$
єэшЄрЁэю ¤ътштрыхэЄхэ юяхЁрЄюЁє т $L_{2}(\T^{d};\mathbb{C}^{\nu})$
єьэюцхэш  эр эхъюЄюЁє■ ьрЄЁшЎє-ЇєэъЎш■ $h(k)$, $k\in\T^{d}$.
╧єёЄ№, ъръ ш т√°х, $\{E_{s}(k)\}_{s=1}^{\nu}$ --- ёюсёЄтхээ√х
чэрўхэш  ьрЄЁшЎ√ $h(k)$, чрэєьхЁютрээ√х ё єўхЄюь ъЁрЄэюёЄш т
яюЁ фъх эхєс√трэш . ╧єёЄ№ $\{\varphi_{s}(k)\}_{s=1}^{\nu}$
--- юЁЄюэюЁьшЁютрээ√щ срчшё шч ёюсёЄтхээ√ї тхъЄюЁют ьрЄЁшЎ√
$h(k)$, юЄтхўр■∙шї ёюсёЄтхээ√ь ўшёырь $\{E_{s}(k)\}_{s=1}^{\nu}$;
$P_{s}(k):=(\cdot,\varphi_{s}(k))\varphi_{s}(k)$, $s=1,\dots,\nu$.
╚ч яЁхфёЄртыхэш  (\ref{KorSlo.1.5}) ёыхфєхЄ ЁртхэёЄтю
\begin{equation}\label{KorSlo.3.2}
\V^{1/2}\mathbb{E}_{+}(\lambda-\mathbb{H})^{-1}\mathbb{E}_{+}\V^{1/2}=\sum_{s=1}^{N-1}\Upsilon_{s}^{*}(\lambda)\Upsilon_{s}(\lambda),\
\ \lambda\in(\Lambda_{+},\Lambda_{-}),
\end{equation}
уфх юяхЁрЄюЁ√
$\Upsilon_{s}(\lambda):\ell_{2}(\zd;\mathbb{C}^{\nu})\to
L_{2}(\T^{d};\mathbb{C}^{\nu})$ юяЁхфхыхэ√ ЁртхэёЄтрьш
\begin{equation*}
\Upsilon_{s}(\lambda):=[(\lambda-E_{s}(k))^{-1/2}P_{s}(k)]\Phi\V^{1/2},\
\ s=1,\dots,N-1,\ \ \lambda\in[\Lambda_{+},\Lambda_{-}).
\end{equation*}
╟фхё№ $\Phi$ --- фшёъЁхЄэюх яЁхюсЁрчютрэшх ╘єЁ№х, юяЁхфхыхээюх
ЁртхэёЄтюь (\ref{KorSlo.1.4}).
\par\noindent
2) ┼ёыш юяхЁрЄюЁ√ $\Upsilon_{s}(\Lambda_{+})$, $s=1,\dots,N-1$
ъюьяръЄэ√, Єю т ёшыє ёююЄэю°хэшщ
$$0\le(\lambda-E_{s}(k))^{-1}\le (\Lambda_{+}-E_{s}(k))^{-1},\ \ s=1,\dots,N-1,\ \ k\in\T^{d},$$
ёяЁртхфышт√ юяхЁрЄюЁэ√х эхЁртхэёЄтр
\begin{equation}\label{KorSlo.3.3}
\Upsilon_{s}^{*}(\lambda)\Upsilon_{s}(\lambda)\le
\Upsilon_{s}^{*}(\Lambda_{+})\Upsilon_{s}(\Lambda_{+}),\ \
s=1,\dots,N-1,\ \ \lambda\in(\Lambda_{+},\Lambda_{-}).
\end{equation}
╚ч (\ref{KorSlo.3.1}), (\ref{KorSlo.3.2}) ш (\ref{KorSlo.3.3})
т√ЄхърхЄ юЎхэър
\begin{equation*}
N_{+}(\lambda,\tau)\le
n(\tau^{-1},\sum_{s=1}^{N-1}\Upsilon_{s}^{*}(\Lambda_{+})\Upsilon_{s}(\Lambda_{+})),\
\ \lambda\in(\Lambda_{+},\Lambda_{-}),
\end{equation*}
ъюЄюЁр , яЁш яхЁхїюфх ъ яЁхфхыє яЁш $\lambda\to\Lambda_{+}+0$,
фрхЄ эхЁртхэёЄтю
\begin{equation}\label{KorSlo.3.4}
N_{+}(\Lambda_{+},\tau)\le
n(\tau^{-1},\sum_{s=1}^{N-1}\Upsilon_{s}^{*}(\Lambda_{+})\Upsilon_{s}(\Lambda_{+})).
\end{equation}
\noindent 3) ╬ЄьхЄшь, ўЄю юяхЁрЄюЁ
$\mathbb{V}=\mathbb{U}V\mathbb{U}^{*}$ хёЄ№ юяхЁрЄюЁ єьэюцхэш  эр
ьрЄЁшЎє-ЇєэъЎш■ т $\zd$
\begin{equation*}
\mathbb{V}(n)=\left(
\begin{matrix}
V(x_{1}+n)&0&\ldots&0\\
0&V(x_{2}+n)&\ldots&0\\
\vdots&\vdots&\ddots&\vdots\\
0&0&\ldots&V(x_{\nu}+n)
\end{matrix}
\right),\ \ n\in\zd.
\end{equation*}
╤ыхфютрЄхы№эю, шч єёыютш  (\ref{KorSlo.1.10}) т√Єхър■Є ёююЄэю°хэш 
\begin{equation}\label{KorSlo.3.5}
\V\in\ell_{p}(\zd),\ \ \|\V\|_{\ell_{p}}^{p}=\|V\|_{\ell_{p}}^{p};
\end{equation}
шч єёыютш  (\ref{KorSlo.1.11}) --- ёююЄэю°хэш 
\begin{equation}\label{KorSlo.3.6}
\V\in\ell_{p,\infty}(\zd),\ \ \|\V\|_{\ell_{p,\infty}}^{p}\le
C(\nu)\|V\|_{\ell_{p,\infty}}^{p};
\end{equation}
эръюэхЎ, шч єёыютш  (\ref{KorSlo.1.12}) --- рёшьЄюЄшър
\begin{equation}\label{KorSlo.3.7}
\mathbb{V}(n)=|n|^{-d/p}\left(\vartheta(\frac{n}{|n|})\mathbbm{1}+o(1)\right),\
\ |n|\to+\infty.
\end{equation}
{\pr \bf ─юърчрЄхы№ёЄтю ЄхюЁхь√ \ref{KorSloTheoremRLC}.} ┬
ёююЄтхЄёЄтшш ё чрьхўрэш ьш яЁхф√фє∙хую яєэъЄр, т єёыютш ї ЄхюЁхь√
\ref{KorSloTheoremRLC} ёяЁртхфышт√ єЄтхЁцфхэш 
(\ref{KorSlo.1.17}$+$) ш (\ref{KorSlo.3.5}). ╤ыхфютрЄхы№эю, ёюуырёэю
яЁхфыюцхэш■ \ref{KorSloPropositionDiscreteCwikelEstimate}, юяхЁрЄюЁ√
$\Upsilon_{s}(\Lambda_{+})\in \SS_{2p,\infty}^{0}$, $s=1,\dots,N-1$,
 ш тхЁэ√ эхЁртхэёЄтр
\begin{equation}\label{KorSlo.3.8}
\|\Upsilon_{s}(\Lambda_{+})\|_{\SS_{2p,\infty}}^{2p}\le
C(p,d,\nu)\|(\Lambda_{+}-E_{s}(\cdot))^{-1}\|_{L_{p,\infty}}^{p}\|V\|_{\ell_{p}}^{p},\
\ s=1,\dots,N-1.
\end{equation}
╬Ўхэъш (\ref{KorSlo.1.18}$+$) ш (\ref{KorSlo.1.19}$+$) т√Єхър■Є шч
(\ref{KorSlo.3.4}), тъы■ўхэшщ
$\Upsilon_{s}(\Lambda_{+})\in\SS_{2p,\infty}^{0}$ ш
(\ref{KorSlo.3.8}).%
\hfill$\square$ \pr {\bf ─юърчрЄхы№ёЄтю ЄхюЁхь√
\ref{KorSloTheoremWRLC}.} ╤юуырёэю чрьхўрэш ь яєэъЄр 1, т єёыютш ї
ЄхюЁхь√ \ref{KorSloTheoremWRLC} ёяЁртхфышт√ єЄтхЁцфхэш 
(\ref{KorSlo.1.16}$+$) ш (\ref{KorSlo.3.6}). ╤ыхфютрЄхы№эю, т ёшыє
яЁхфыюцхэш  \ref{KorSloPropositionCwikelTypeEstimate}, тхЁэ√
ёююЄэю°хэш 
\begin{equation}\label{KorSlo.3.9}
\Upsilon_{s}(\Lambda_{+})\in\SS_{2p,\infty},\ \
\|\Upsilon_{s}(\Lambda_{+})\|_{\SS_{2p,\infty}}^{2p}\le
C(p,d,\nu,\varkappa)\|(\Lambda_{+}-E_{s}(\cdot))^{-1}\|_{L_{\varkappa}}^{p}\|V\|_{\ell_{p,\infty}}^{p},\
\
\end{equation}
яЁш тёхї $ s=1,\dots,N-1$. ╬Ўхэър (\ref{KorSlo.1.20}$+$) т√ЄхърхЄ шч
(\ref{KorSlo.3.4}) ш (\ref{KorSlo.3.9}).

┼ёыш фюяюыэшЄхы№эю $V\in\ell_{p,\infty}^{\,0}(X)$, Єю
$\V\in\ell_{p,\infty}^{\,0}(\zd)$, р яюЄюьє (т ёшыє яЁхфыюцхэш 
\ref{KorSloPropositionCwikelTypeEstimate}) ёяЁртхфышт√ тъы■ўхэш 
$\Upsilon_{s}(\Lambda_{+})\in\SS_{2p,\infty}^{0}$,
$s=1,\dots,N-1$, р чэрўшЄ ш тъы■ўхэшх
$\sum_{s=1}^{N-1}\Upsilon_{s}^{*}(\Lambda_{+})\Upsilon_{s}(\Lambda_{+})\in\SS_{p,\infty}^{0}$.
╧юёыхфэхх тьхёЄх ё (\ref{KorSlo.3.4}) ш ючэрўрхЄ
$N_{+}(\Lambda_{+},\tau)=o(\tau^{p})$, $\tau\to+\infty$.
\hfill$\square$%
{\pr \bf ─юърчрЄхы№ёЄтю ЄхюЁхь√ \ref{KorSloTheoremAsymptotic}.
╤ыєўрщ
$\lambda\in(\Lambda_{+},\Lambda_{-})$.} 
╩ръ ш т√°х, ЁрчсхЁхь сюыхх ЄЁєфэ√щ ёыєўрщ {\it тэєЄЁхээхщ\/}
ыръєэ√ $(\Lambda_{+},\Lambda_{-})$ т ёяхъЄЁх юяхЁрЄюЁр $H$,
єфютыхЄтюЁ ■∙хщ яЁш эхъюЄюЁюь $N\in\{2,\dots,\nu\}$ єёыютш■
(\ref{KorSlo.1.7}). ╧Ёш Єръшї єёыютш ї ЁртхэёЄтю
(\ref{KorSlo.1.15}$+$) яЁшэшьрхЄ тшф
\begin{equation*}
\Gamma_{p}^{+}(\lambda):=
(2\pi)^{-d}d^{-1}\sum_{s=1}^{N-1}\int_{\T^{d}}\left(\lambda-E_{s}(k)\right)^{-p}dk
\int_{\mathbb{S}^{d-1}}\vartheta^{p}(\theta)dS(\theta),\ \
\lambda\in[\Lambda_{+},\Lambda_{-}].
\end{equation*}
╚ч ЁртхэёЄтр (\ref{KorSlo.1.14}$+$) ш трЁшрЎшюээюую ётющёЄтр
(\ref{KorSlo.2.1}) т√Єхър■Є юЎхэър ётхЁїє (\ref{KorSlo.3.1}) ш
юЎхэър {\it ёэшчє\/}
\begin{multline}\label{KorSlo.3.10}
N_{+}(\lambda,\tau)\ge n_{+}(\tau^{-1},E_{+}V^{1/2}(\lambda
I-H)^{-1}V^{1/2}E_{+})=\\=n_{+}(\tau^{-1},\mathbb{E}_{+}\V^{1/2}(\lambda
I-\mathbb{H})^{-1}\V^{1/2}\mathbb{E}_{+}),\ \
\lambda\in(\Lambda_{+},\Lambda_{-}),\ \ \tau>0.
\end{multline}
┬ єёыютш ї ЄхюЁхь√ \ref{KorSloTheoremAsymptotic} юяхЁрЄюЁ
$\V^{1/2}$ яЁшэрфыхцшЄ ъырёёє $\SS_{2p,\infty}$. ╥ръшь юсЁрчюь, шч
(\ref{KorSlo.3.1}) ш (\ref{KorSlo.3.10}) ёыхфєхЄ, ўЄю фы 
юсюёэютрэш  (\ref{KorSlo.1.21}$+$) фюёЄрЄюўэю яЁютхЁшЄ№ ЁртхэёЄтр
\begin{equation}\label{KorSlo.3.11}
\mathfrak{d}_{p}^{+}(\mathbb{E}_{+}\mathbb{V}^{1/2}(\lambda
I-\mathbb{H})^{-1}\mathbb{V}^{1/2}\mathbb{E}_{+})=\mathfrak{D}_{p}^{+}(\mathbb{V}^{1/2}\mathbb{E}_{+}(\lambda
I-\mathbb{H})^{-1}\mathbb{E}_{+}\mathbb{V}^{1/2})=\Gamma_{p}^{+}(\lambda),\
\ \lambda\in(\Lambda_{+},\Lambda_{-}).
\end{equation}
╤юуырёэю (\ref{KorSlo.1.5}) шьххЄ ьхёЄю яЁхфёЄртыхэшх
$\mathbb{E}_{+}=\Phi_{}^{*}[P_{+}(k)]\Phi_{}$, уфх
$P_{+}(k):=\sum_{s=1}^{N-1}P_{s}(k)$; яЁш ¤Єюь т єёыютш ї ЄхюЁхь√
\ref{KorSloTheoremAsymptotic} ёяЁртхфыштр рёшьяЄюЄшър
(\ref{KorSlo.3.7}).
╚ч яЁхфыюцхэш  \ref{KorSloPropositionApproxComm} т√ЄхърхЄ
тъы■ўхэшх
$$\mathbb{E}_{+}\mathbb{V}^{1/2}-\V^{1/2}\mathbb{E}_{+}\in\SS_{2p,\infty}^{0},$$
р яюЄюьє тхЁэю
\begin{equation}\label{KorSlo.3.12}
\mathbb{E}_{+}\mathbb{V}^{1/2}(\lambda
I-\mathbb{H})^{-1}\mathbb{V}^{1/2}\mathbb{E}_{+}-\mathbb{V}^{1/2}\mathbb{E}_{+}(\lambda
I-\mathbb{H})^{-1}\mathbb{E}_{+}\mathbb{V}^{1/2}\in\SS_{p,\infty}^{0},\
\ \lambda\in(\Lambda_{+},\Lambda_{-}).
\end{equation}
╤ єўхЄюь (\ref{KorSlo.3.12}) ьюцэю чрьхэшЄ№ (\ref{KorSlo.3.11}) эр
¤ътштрыхэЄэюх т√Ёрцхэшх
\begin{equation}
\label{KorSlo.3.13}
\mathfrak{d}_{p}(\mathbb{V}^{1/2}\mathbb{E}_{+}(\lambda
I-\mathbb{H})^{-1}\mathbb{E}_{+}\mathbb{V}^{1/2})=\mathfrak{D}_{p}(\mathbb{V}^{1/2}\mathbb{E}_{+}(\lambda
I-\mathbb{H})^{-1}\mathbb{E}_{+}\mathbb{V}^{1/2})=\Gamma_{p}^{+}(\lambda),\
\ \lambda\in(\Lambda_{+},\Lambda_{-}).
\end{equation}
╧юёъюы№ъє юяхЁрЄюЁ $\mathbb{E}_{+}(\lambda
I-\mathbb{H})^{-1}\mathbb{E}_{+}$ яюыюцшЄхы№эю юяЁхфхыхэ,
ёююЄэю°хэш  (\ref{KorSlo.3.13}) Ёртэюёшы№э√ ЁртхэёЄтрь
\begin{multline}\label{KorSlo.3.14}
\mathfrak{d}_{p}(\mathbb{E}_{+}(\lambda
I-\mathbb{H})^{-1/2}\mathbb{E}_{+}\mathbb{V}\mathbb{E}_{+}(\lambda
I-\mathbb{H})^{-1/2}\mathbb{E}_{+})=\\
=\mathfrak{D}_{p}(\mathbb{E}_{+}(\lambda
I-\mathbb{H})^{-1/2}\mathbb{E}_{+}\mathbb{V}\mathbb{E}_{+}(\lambda
I-\mathbb{H})^{-1/2}\mathbb{E}_{+})=\Gamma_{p}^{+}(\lambda),\ \
\lambda\in(\Lambda_{+},\Lambda_{-}).
\end{multline}
┬ ёшыє (\ref{KorSlo.1.5}), ёююЄэю°хэш  (\ref{KorSlo.3.14}) ьюцэю
чряшёрЄ№ т тшфх
\begin{multline}\label{KorSlo.3.15}
\mathfrak{d}_{p}([P_{+}(k)(\lambda\mathbf{1}-h(k))^{-1/2}P_{+}(k)]\Phi_{}\mathbb{V}
\Phi_{}^{*}[P_{+}(k)(\lambda\mathbf{1}-h(k))^{-1/2}P_{+}(k)])=\\
=\mathfrak{D}_{p}([P_{+}(k)(\lambda\mathbf{1}-h(k))^{-1/2}P_{+}(k)]
\Phi_{}\mathbb{V}\Phi_{}^{*}[P_{+}(k)(\lambda\mathbf{1}-h(k))^{-1/2}P_{+}(k)])=\Gamma_{p}^{+}(\lambda),\
\
\end{multline}
яЁш $\lambda\in(\Lambda_{+},\Lambda_{-})$. ╬ёЄрхЄё  чрьхЄшЄ№, ўЄю
(\ref{KorSlo.3.15}) т√ЄхърхЄ эхяюёЁхфёЄтхээю шч (\ref{KorSlo.2.15}).

╤ыєўрщ яюыєсхёъюэхўэющ ыръєэ√ ЁрчсшЁрхЄё  рэрыюушўэю (Єюы№ъю
эхёъюы№ъю яЁю∙х).
\par{\bf ╤ыєўрщ $\lambda=\Lambda_{+}$.} ╧ЁютхЁшь ЇюЁьєыє (\ref{KorSlo.1.21}$+$) эр
ыхтюь ъЁр■ ыръєэ√ яЁш єёыютшш (\ref{KorSlo.1.16}$+$). ╩ръ ш т√°х,
ЁрчсхЁхь ёыєўрщ тэєЄЁхээхщ ыръєэ√ $(\Lambda_{+},\Lambda_{-})$ т
ёяхъЄЁх юяхЁрЄюЁр $H$, єфютыхЄтюЁ ■∙хщ яЁш эхъюЄюЁюь
$N\in\{2,\dots,\nu\}$ єёыютш■ (\ref{KorSlo.1.7}). ╚ч
(\ref{KorSlo.1.16}$+$) т√ЄхърхЄ ёююЄэю°хэшх
$\Gamma_{p}^{+}(\lambda)\to\Gamma_{p}^{+}(\Lambda_{+})$,
$\lambda\to\Lambda_{+}+0$. ╤ фЁєующ ёЄюЁюэ√, яЁш єёыютшш
(\ref{KorSlo.1.16}$+$) фы  ы■сюую $s=1,\dots,N-1$ ёяЁртхфышт√
єЄтхЁцфхэш 
\begin{equation}\label{KorSlo.3.16}
\begin{gathered}
\|(\lambda-E_{s}(\cdot))^{-1/2}P_{s}(\cdot)\|_{2\varkappa}\le\|(\Lambda_{+}-E_{s}(\cdot))^{-1/2}P_{s}(\cdot)\|_{2\varkappa},\
\ \lambda\in(\Lambda_{+},\Lambda_{-}),\\
(\lambda-E_{s}(\cdot))^{-1/2}P_{s}(\cdot)\to(\Lambda_{+}-E_{s}(\cdot))^{-1/2}P_{s}(\cdot),\
\ \lambda\to\Lambda_{+}+0,\ \ \text{т}\ \ L_{2\varkappa}(\T^{d}).
\end{gathered}
\end{equation}
╤юуырёэю яЁхфыюцхэш■ \ref{KorSloPropositionCwikelTypeEstimate}, шч
(\ref{KorSlo.3.16}) ёыхфєхЄ, ўЄю юяхЁрЄюЁ $\Upsilon_{s}(\lambda)$
ёїюфшЄё  яЁш $\lambda\to\Lambda_{+}+0$ т ъырёёх $\SS_{2p,\infty}$
ъ юяхЁрЄюЁє $\Upsilon_{s}(\Lambda_{+})$, $s=1,\dots,N-1$. ╧Ёш ¤Єюь
ътрчшэюЁьр $\|\Upsilon_{s}(\lambda)\|_{2\varkappa,\infty}$,
$\lambda\in(\Lambda_{+},\Lambda_{-})$, ЁртэюьхЁэю юуЁрэшўхэр.
╤ыхфютрЄхы№эю (ёь. яЁхфыюцхэшх \ref{KorSloPropositionWOGelder} ш
(\ref{KorSlo.3.2})), юяхЁрЄюЁ
$\mathbb{V}^{1/2}\mathbb{E}_{+}(\lambda
I-\mathbb{H})^{-1}\mathbb{E}_{+}\mathbb{V}^{1/2}$ ёїюфшЄё  яЁш
$\lambda\to\Lambda_{+}+0$ т ъырёёх $\SS_{p,\infty}$ ъ юяхЁрЄюЁє
$\sum_{s=1}^{N-1}\Upsilon_{s}^{*}(\Lambda_{+})\Upsilon_{s}(\Lambda_{+})$,
р яюЄюьє
\begin{equation}\label{KorSlo.3.17}
\mathfrak{d}_{p}\left(\sum_{s=1}^{N-1}\Upsilon_{s}^{*}(\Lambda_{+})\Upsilon_{s}(\Lambda_{+})\right)=\mathfrak{D}_{p}
\left(\sum_{s=1}^{N-1}\Upsilon_{s}^{*}(\Lambda_{+})\Upsilon_{s}(\Lambda_{+})\right)=\Gamma_{p}^{+}(\Lambda_{+}).
\end{equation}
╩ръ ш т√°х, т єёыютш ї ЄхюЁхь√ \ref{KorSloTheoremAsymptotic}
ёяЁртхфыштю эхЁртхэёЄтю (\ref{KorSlo.3.4}). ╤ыхфютрЄхы№эю (т ёшыє
ьюэюЄюээюёЄш $N_{+}(\cdot,\tau)$), ёяЁртхфышт√ эхЁртхэёЄтр
\begin{equation}\label{KorSlo.3.18}
N_{+}(\lambda,\tau)\le N_{+}(\Lambda_{+},\tau)\le
n(\tau^{-1},\sum_{s=1}^{N-1}\Upsilon_{s}^{*}(\Lambda_{+})\Upsilon_{s}(\Lambda_{+})),\
\ \lambda\in(\Lambda_{+},\Lambda_{-}),\ \ \tau>0.
\end{equation}
╚ч (\ref{KorSlo.3.17}) ш (\ref{KorSlo.3.18}) т√Єхър■Є ёююЄэю°хэш 
\begin{equation}\label{KorSlo.3.19}
\Gamma_{p}^{+}(\lambda)\le\liminf_{\tau\to+\infty}\tau^{-p}N_{+}(\Lambda_{+},\tau)\le\limsup_{\tau\to+\infty}\tau^{-p}N_{+}(\Lambda_{+},\tau)\le\Gamma_{p}^{+}(\Lambda_{+}),\
\ \lambda\in(\Lambda_{+},\Lambda_{-}).
\end{equation}
╧хЁхїюф  т (\ref{KorSlo.3.19}) ъ яЁхфхыє яЁш
$\lambda\to\Lambda_{+}+0$, яюыєўшь
\begin{equation*}
\lim_{\tau\to+\infty}\tau^{-p}N_{+}(\Lambda_{+},\tau)=\Gamma_{p}^{+}(\Lambda_{+}).
\end{equation*}
\hfill$\square$

\end{document}